\documentclass[11pt]{amsart}
\usepackage{geometry}                
\usepackage{graphicx}
\usepackage{amssymb}
\usepackage{epstopdf}
\usepackage{mathrsfs}
\usepackage{hyperref}

\newtheorem{theorem}{Theorem}[section]
\newtheorem{lemma}[theorem]{Lemma}

\newtheorem{corollary}[theorem]{Corollary}
\theoremstyle{definition}
\newtheorem{definition}[theorem]{Definition}
\newtheorem{example}[theorem]{Example}

\newtheorem{remark}[theorem]{Remark}

\newcommand{\n}{\|} 

\newcommand{\la}{\langle}
\newcommand{\ra}{\rangle}
\newcommand{\LL}{\mathcal{L}}
\usepackage{cancel}
\usepackage{xcolor}

\title[Smoothing and Strichartz estimates for degenerate equations]{Smoothing and Strichartz estimates for degenerate Schr\"{o}dinger-type equations}
\author[S. Federico]{Serena Federico}
\thanks{This project has received funding from the European Union’s Horizon 2020 research and innovation programme under the Marie Sk\l odowska-Curie grant agreement No 838661.}
\address{
Serena Federico
\endgraf
Department of Mathematics: Analysis, Logic and Discrete Mathematics,
\endgraf
Ghent University,
Belgium}\email[Serena Federico ]{serena.federico@ugent.be}
\author[M. Ruzhansky]{Michael Ruzhansky}
\thanks{MR was supported by the FWO Odysseus 1 grant G.0H94.18N: Analysis and Partial Differential Equations and in parts by the EPSRC Grant EP/R003025/1 and by the Leverhulme Research Grant RPG-2017-151.}
\address{
Michael Ruzhansky
\endgraf
Department of Mathematics: Analysis, Logic and Discrete Mathematics,
\endgraf
Ghent University,
Belgium
\endgraf
and
\endgraf
School of Mathematical Sciences\endgraf
Queen Mary University of London\endgraf
United Kingdom}\email[Michael Ruzhansky]{michael.ruzhansky@ugent.be}
\keywords{Smoothing effect; Strichartz estimates; Schr\"{o}dinger operators; Comparison principles;}

\begin{document}
\maketitle

\begin{abstract}
In this paper we focus on the validity of some fundamental estimates for time-degenerate Schr\"{o}dinger-type operators.  On one hand we derive global homogeneous smoothing estimates for operators of any order by means of suitable comparison principles (that we shall obtain here). 
On the other hand, we prove weighted Strichartz-type estimates for time-degenerate Schr\"{o}dinger operators and apply them to the local well-posedness of the semilinear Cauchy problem. Most of our results apply to nondegenerate operators as well, recovering, in these cases, the well-known standard results.
\end{abstract}
\tableofcontents

\section{Introduction}
In what follows we shall discuss the global homogeneous smoothing effect for some time-degenerate Schr\"{o}dinger-type operators. Moreover, for the same class of operators, we shall derive some weighted Strichartz estimates which will be employed to prove local well-posedness results for the associated semilinear Cauchy problem.

We are interested in operators of the form

\begin{equation}
\label{A}
\mathcal{L}_{b,a}=\partial_t -iB(t) a(D),
\end{equation}
where $a(D)$ is a Fourier multiplier of order $m$ with symbol denoted by $a(\xi)$,
 $B\in C(\mathbb{R})$ and $b=b(t)=\int_0^t B(s)ds$. 
 As usual, here $a(D)f=\mathcal{F}^{-1}(a(\xi)\mathcal{F}f)$, where $\mathcal{F}$ and $\mathcal{F}^{-1}$ are the Fourier and inverse Fourier transforms, respectively. 
 We are interested in the case when $B$ is such that $B(0)=0$, that is the time-degenerate case, but some of our results  will be applicable to more general situations, that is, for instance, when $B$ is vanishing at finitely many points and, also, when $B$ vanishes at infinitely many points. However, also the case when $B$ is non-vanishing is covered by our results, and, in the particular, when $B(t)=1$, that is when $b(t)=t$, one can recover classical results for constant coefficients Schr\"{o}dinger equations.
The reason for the notation $\LL_{b,a}$ is dictated by the fact that by using Fourier analysis one immediately has that the solution at time $t$ of the homogeneous Cauchy problem 
\begin{equation}
\label{IVP}
\left\{ \begin{array}{l}
 \LL_{b,a} u=0, \\
u(s,x)=u_s(x),\\
\end{array}
\right .
\end{equation}
is of the form
$$u(t,x)=e^{i (b(t)-b(s)) a(D)} u_s(x):= \int_{\mathbb{R}^n}e^{ix\cdot \xi +i (b(t)-b(s))a(\xi)}\widehat{u_s}(\xi)d\xi,$$
where $\widehat{u_s}$ stands for the (space-)Fourier transform of $u_s$.

The problems we address in this paper have been intensively studied in the nondegenerate setting. The homogeneous smoothing effect, which describes a gain of smoothness of the homogeneous solution of the IVP with respect to the smoothness of the initial data, and the inhomogeneous smoothing effect, describing a gain of smoothness of the solution of the inhomogeneous IVP with respect to the regularity of the inhomogeneous data, have been shown for general dispersive equations. The smoothing effect for the Korteweg-de Vries equation was first shown by Kato in \cite{K}, while, for Schr\"{o}dinger equations, the homogeneous smoothing effect was simultaneously established by Sj\"{o}lin \cite{Sj} and Vega \cite{V}.  The homogeneous smoothing effect for general dispersive equations was proved by Constantin and Saut in \cite{CS}, and, later, improved and generalized by Kenig, Ponce and Vega in \cite{KPV2}, where also the inhomogeneous smoothing effect for Schr\"{o}dinger equations was shown, and also by  Ben-Artzi and Klainerman in \cite{BAK}.
Other global smoothing estimates were also derived by Chihara \cite{C}, Kato and Yajima \cite{Kato-Yajima}, Linares and Ponce \cite{LP}, Sugimoto \cite{S}, Walther \cite{W,W2} and many others.

The results mentioned above deal with {\it constant coefficients} operators. In the variable coefficients case the results obtained so far mainly concern with the nondegenerate setting.
The smoothing effect for generalized and ultrahyperbolic variable coefficients Schr\"{o}dinger operators was proved by Kenig et al.  in \cite{KPV4} and \cite{KPRV} respectively. In these papers the local smoothing property was also employed to prove local well-posedness results for the associated nonlinear Cauchy problem.  Important contributions are due to Doi (see \cite{D}), who also proved local smoothing results in the general manifold setting. For variable coefficients nondegenerate Schr\"{o}dinger operators smoothing results were recently proved by Marzuola, Metcalfe and Tataru in \cite{MMT} (see also reference therein).

In the context of degenerate operators, few results have been obtained about the smoothing effect. Staffilani and the first author derived in \cite{FS} the local smoothing effect for time-degenerate Schr\"{o}dinger operators of the form considered here, and used them to obtain the local well-posedness of the associated nonlinear IVP. It is also worth to mention that degenerate operators of the form considered in \cite{FS} were previously analized by Cicognani and Reissig in \cite{CR}, who studied the local well-posedness of the linear (homogeneous) IVP both in Sobolev and Gevrey spaces.

On the other side, that is on the side of Strichartz-type estimates, there is a lot of literature about these estimates for constant and variable coefficients partial differential operators. 
These estimates, as opposite to the smoothing ones, show a gain of integrability instead of a gain of smoothness of the solution of the IVP.
The starting point in this direction was the original paper by Strichartz \cite{STR}, later on extended by Ginibre and Velo in \cite{GV}. The resolution of the endpoint case, instead, is due to Keel and Tao \cite{KT}.
After that, many results have been obtained, once again concerning the nondegenerate cases, both in the Euclidean and in the manifold setting. We will only mention some of them here. For instance, in the variable coefficients case in the Euclidean setting, one has results in \cite{ST, RZ, MMT, M,EGS,EGS2} (see also references therein), whereas, in the manifold setting, one has results in \cite{B,BGT,RT,HTW} (see also references therein).

Our first goal here is to investigate the homogeneous global smoothing effect of the family of operators 
\begin{equation*}
\{e^{i b(t) a(D)}\}_{t\in\mathbb{R}}
\end{equation*}
which corresponds to the one-parameter family of unitary operators giving the solution of \eqref{IVP} in the case $s=0$. 
Notice that, depending on the properties of the function $b$, 
we have that the family of unitary operators described above is a group or a semigroup.
The presence of the function $B$ in \eqref{A} makes the operator degenerate when the function vanishes at some point of the time domain. However, nondegenerate cases are covered by our results as well.
Note also that, in the case $b(t)=t$, one has classical smoothing results for the solutions of the Cauchy problem for general dispersive equations of the form \eqref{A}, while, in the case $b(t)=t^\alpha, \alpha>0$, and $m=2$, some weighted local smoothing estimates have been recently derived in \cite{FS}. We remark that, with respect to the homogeneous smoothing effect proved in \cite{FS} where a finite time domain is considered, here the homogeneous smoothing effect holds globally in time.

We stress that weighted smoothing estimates hold for non dispersive equations as well, and that our class covers such cases too. For non dispersive equations non standard smoothing results were obtained by Sugimoto and the second author in \cite{RS2}, where it was conjectured that, given an operator of the form $L=\partial_t+ia(t,D_x)$, an invariant estimate of the form
\begin{equation}
\label{conj}
\n \la x\ra^{-s}|\nabla_x a(t,D_x)|^{1/2}e^{i\int_0^ta(s,D_x)ds}\varphi\n_{L^2(\mathbb{R}_t\times\mathbb{R}^n_x)}\leq C\n \varphi\n_{L^2(\mathbb{R}^n)},
\end{equation}
holds. In fact, following this idea, they proved a weighted homogeneous smoothing effect for time dependent equations of the form \eqref{A} with $B(t)> 0$. 
The smoothing result we prove here extends the latter to the more general degenerate setting, also supporting the above conjecture. Estimates of the form \eqref{conj} were called universal estimates in \cite{RS2}.

The homogeneous smoothing effect for $\LL_{b,|D_x|^m}$ can be derived from the standard case after application of Lemma \ref{T1} below. In fact, the strength of the Lemma \ref{T1} is to reduce the possibly degenerate case $\LL_{b,|D_x|^m}$ to the nondegenerate case $\LL_{t,|D_x|^m}$.  
However, we obtain the homogeneous smoothing effect for $\LL_{b,|D_x|^m}$  as a consequence of the validity of comparison principles between operators of the general form $\LL_{b,a}$. 

Comparison principles for operators in our class are the subject of Section \ref{CP}. In particular, Section \ref{CP}  will be devoted to a suitable refinement of the comparison principles introduced by Sugimoto and the second author in \cite{RS} (more abstract spectral comparison principles have been derived in \cite{RSBA}).

 In Section \ref{ACP} we will then show the application of the tools developed in Section \ref{CP} to the derivation of smoothing estimates for the operators $\LL_{b,a}$ with $a(D)=|D_x|^m$ or $a(D)=D_{x_j}|D_{x'}|^{m-1}$. We will show that the result can be derived by comparison with the standard case. The global homogeneous smoothing effect, which is given by means of a weight related to the function $b$, is contained in Theorem \ref{Thm2} of Section \ref{ACP}. In Section \ref{ACP} we also state in Theorem \ref{GSE} the natural extension to the degenerate setting of the classical homogeneous smoothing effect for Schr\"{o}dinger operators ($m=2$).
In addition, once again by comparison with the classical case, we derive global weighted smoothing estimates for operators of the form $\LL_{b,|D_x|^m}$ equivalent to those obtained by Chihara in \cite{C}, by Sugimoto in \cite{S} and by Walther in \cite{W} for operators of the form $\LL_{t,|D_x|^m}$.

Weighted Strichartz-type estimates for time-degenerate Schr\"{o}dinger operators $\LL_{b,\Delta}$ are studied in Section \ref{Strichartz}. We first derive global weighted Strichartz estimates with the standard range of admissible exponents $(q,p)$. Then we obtain (different) local Strichartz estimates in the nonendpoint case. 

Finally, in Section \ref{LWP}, we employ the local in time Strichartz estimates of Section \ref{Strichartz} to prove the local well-posedness of the semilinear IVP with a nonlinearity of the form $N=N(B,u)=|B||u|^{p-1}u$. 
We will conclude Section \ref{LWP} by giving some examples of operators to which our smoothing and local well-posedness results apply.

\vspace{0.3cm}

\noindent{\bf{Notations}.} In what follows, to simplify the notation, we shall write $b'(t)$ in place of $B(t)$ in \eqref{A}, and, sometimes, we will write equation \eqref{A} in the form
$$\mathcal{L}_{b,a}=\partial_t -ib'(t) a(D),$$
where $b\in C^1(\mathbb{R})$ will be, according to the previous definition, such that $b(0)=0$.
Moreover we shall write $\lim_{t\rightarrow \infty}b(t)=\infty$ when either $\lim_{t\rightarrow \pm\infty}b(t)=\pm\infty$ or $\lim_{t\rightarrow \pm\infty}b(t)=\mp\infty$.

\section{Comparison principles}\label{CP}
In this section we shall derive comparison principles for time-degenerate equations. 
In particular, we shall make use of the key result proved  in \cite{RS} to obtain a weighted version of comparison principles suitable to our case.

Given two partial differential operators $a(D)$ and $\tilde{a}(D)$ (not necessarily elliptic and eventually $a=\tilde{a}$) with symbols $a(\xi)$ and $\tilde{a}(\xi)$ respectively, and two functions $b, f\in C^1(\mathbb{R})$ vanishing at 0 (possibly $b=f$), we shall denote by $\LL_{b,a}$, and $\LL_{f,\tilde{a}}$  the operators
$$\mathcal{L}_{b,a}=\partial_t -ib'(t) a(D) \quad\text{and}\quad \mathcal{L}_{f, \tilde{a}}=\partial_t -i f'(t)\tilde{a}(D).$$
Then we consider the initial value problems (IVP)
\begin{equation*}
\left\{ \begin{array}{l}
 \LL_{b,a} u=0, \\
u(0,x)=\varphi(x),\\
\end{array}
\right .
\end{equation*}
and
\begin{equation*}
\left\{ \begin{array}{l}
 \LL_{f,\tilde{a}}v=0, \\
v(0,x)=\varphi(x),\\
\end{array}
\right .
\end{equation*}
where $t\in\mathbb{R}$ and $x\in\mathbb{R}^n$, whose solution operators are, respectively, given by
\begin{equation}
W_{b,a}(t)u_0(x)=e^{ib(t) a(D)}\varphi(x):= \int_{\mathbb{R}^n}e^{i b(t) a(\xi)+ix\cdot \xi}\widehat{\varphi}(\xi)d\xi,
\end{equation}
and 
\begin{equation}
W_{f,\tilde{a}}(t)v_0(x)=e^{i f(t)\tilde{a}(D)}\varphi(x):= \int_{\mathbb{R}^n}e^{i f(t) \tilde{a}(\xi)+ix\cdot \xi}\widehat{\varphi}(\xi)d\xi.
\end{equation}

Our goal is to be able to compare smoothing estimates for $\LL_{b,a}$ and $\LL_{f,\tilde{a}}$, so that, by suitably choosing $f$ and $\tilde{a}$, we can derive the smoothing effect for $\LL_{b,a}$ from that of $\LL_{f,\tilde{a}}$.

In order to obtain comparison principles for time-degenerate equations of the form considered above, we will need the following lemma.

\begin{lemma}\label{T1}
Let $n,p\geq 1$,  and let $b\in C^1(\mathbb{R})$ be such that $b(0)=0$. Let also $a\in C^1(\mathbb{R}^n)$ and $\sigma\in C^0(\mathbb{R}^n)$.
Then we have that
\vspace{0,3cm}

\item[(i)] If $b$ is strictly monotone then, for all $u_0\in L^{2}_x(\mathbb{R}^n)$,
\begin{equation*}
\n |b'(t)|^{1/p}\sigma(D_x)e^{i b(t) a(D)} u_0\n_{L^p_t(\mathbb{R})}\leq \n \sigma(D_x)e^{i t a(D)} u_0\n_{L^p_t(\mathbb{R})},\end{equation*}
where equality holds if $b$ also satisfies $\lim_{t\rightarrow \infty}b(t)=\infty$. \\
\item[(ii)] If $b'$ is such that $\sharp\{t\in\mathbb{R}; b'(t)=0\}=k\geq 1$, then, for all $u_0\in L^{2}_x(\mathbb{R}^n)$, 
\begin{equation*}
\n |b'(t)|^{1/p}\sigma(D_x)e^{i b(t) a(D)} u_0\n_{L^p_t(\mathbb{R})}\leq C(k)
\n \sigma(D_x) e^{i t a(D)} u_0\n_{L^p_t(\mathbb{R})}
\end{equation*}
with $C(k)=(k+1)^{1/p}$.\\
\item[(iii)] If $b'$ is such that there exists an increasing sequence of positive times $\{t_k\}_{k\in\mathbb{N}}$ and a decreasing sequence of negative times $\{t'_k\}_{k\in\mathbb{N}}$ such that $b'(t_k)=b'(t'_k)=0$, and we set $t_0=t'_0=0$, with $b'$ possibly vanishing at $0$ as well, then, for any function $c\in C(\mathbb{R})$ such that $\sum_{k=0}^\infty \sup_{[t_k,t_{k+1})}|c(t)|,\, \sum_{k=0}^\infty \sup_{(t'_{k+1},t'_{k}]}|c(t)|\leq C<\infty$, we have, for all $u_0\in L^2_x(\mathbb{R}^n)$,
\begin{equation*}
\n |c(t)b'(t)|^{1/p}\sigma(D_x)e^{i b(t) a(D)} u_0\n_{L^p_t(\mathbb{R})}\leq (2C)^{1/p}
\n \sigma(D_x) e^{i t a(D)} u_0\n_{L^p_t(\mathbb{R})}.
\end{equation*}
\end{lemma}

\proof[Proof of (i)]
By using the change of variables $s=b(t)$ and denoting by $A$ the set $A:=b^{-1}(\mathbb{R})$, we have
$$\n |b'(t)|^{1/p}\sigma(D_x)e^{i b(t) a(D_x)} u_0\n^p_{L^p_t(\mathbb{R})}=
\int_{\mathbb{R}}|b'(t)|\, | \sigma(D_x)e^{i b(t) a(D_x)} u_0|^p dt$$
$$= \int_{A}\frac{|b'(b^{-1}(s))|}{b'(b^{-1}(s))} | \sigma(D_x)e^{i s a(D_x)} u_0|^p ds$$
$$\leq  \int_{\mathbb{R}}|\sigma(D_x)e^{i s a(D_x)} u_0|^p ds=\n \sigma(D_x)e^{i t a(D_x)} u_0\n^p_{L^p_t},$$
where we used that $b$ is strictly monotone. Note that, if $b$ satisfies $\lim_{t\rightarrow \infty}b(t)=\infty$, that is,  $A=(\mp\infty,\pm \infty)$ when $b$ if strictly increasing or strictly decreasing respectively, then we have an equality. 
\endproof

\proof[Proof of (ii)]
In this case we proceed by splitting the time domain into regions where the change of variables in time is allowed.  
Then, denoting by $I_0=(-\infty, t_1]$, $I_k=[t_k,+\infty)$ and by $I_j=[t_j,t_{j+1}]$ for $j=1,...,k-1$, we have
$$\n|b'(t)|^{1/p}\sigma(D_x)e^{i b(t) a(D_x)} u_0\n^p_{L^p_t(\mathbb{R})}= 
\sum_{j=0}^k\int_{I_j} |b'(t)| |\sigma(D_x)e^{i b(t) a(D_x)} u_0|^pdt$$
$$\underset{s=b(t)}{=} \sum_{j=0}^k\int_{I'_j} \frac{|b'(b^{-1}(s))|}{b'(b^{-1}(s))} | \sigma(D_x)e^{i s a(D_x)} u_0|^pds$$
$$\leq (k+1)\int_\mathbb{R} | \sigma(D_x)e^{i s a(D_x)} u_0|^p ds=(k+1)\n \sigma(D_x)e^{i t a(D_x)}u_0\n^p_{L^p_t},$$
which gives (ii).
\endproof

\proof[Proof of (iii)]
We shall adopt here the following notations: $I_j=[t_j,t_{j+1})$, $I'_j=(t'_{j+1},t_j']$ will be positive and negative time intervals respectively, $c_j:=\sup_{I_j}|c|$, and  $c'_j:=\sup_{I'_j}|c|$. By using the properties of the function $c$ we have

$$\n |c(t)b'(t)|^{1/p}\sigma(D_x)e^{i b(t) a(D)} u_0\n^p_{L^p_t(\mathbb{R})}$$ 
$$ =\sum_{j=0}^\infty \int_{I_j} |c(t) b'(t)| | \sigma(D_x)e^{i b(t) a(D_x)} u_0|^pdt+\sum_{j=0}^\infty \int_{I'_j}|c(t) b'(t)| | \sigma(D_x)e^{i b(t) a(D_x)} u_0|^p dt$$
$$\leq \sum_{j=0}^\infty c_j  \int_{I_j} | b'(t)| | \sigma(D_x)e^{i b(t) a(D_x)} u_0|^pdt+ \sum_{j=0}^\infty c'_j  \int_{I'_j} | b'(t)|| \sigma(D_x)e^{i b(t) a(D_x)} u_0|^pdt$$
$$ \underset{b(t)=s}{\leq}\sum_{j=0}^\infty c_j \n\sigma(D_x)e^{i s a(D_x)} u_0\n_{L^p_s}^p+\sum_{j=0}^\infty c'_j \n\sigma(D_x)e^{i s a(D_x)} u_0\n_{L^p_s}^p$$
$$\leq 2C\n \sigma(D_x)e^{i s a(D_x)} u_0\n_{L^p_t}^p,$$
which gives (iii).
\endproof

We now state, in this setting, the theorem from which a series of comparison principles will be derived. This theorem is the suitable generalization of Theorem 2.1 in \cite{RS} in our setting.

\begin{theorem}\label{TC2}
Let $b\in C^1(\mathbb{R)}$ be such that $b(0)=0$ and $a\in C^1(\mathbb{R}^n)$ such that, for almost all $\xi'=(\xi_2,...,\xi_n)\in\mathbb{R}^{n-1}$, $a(\xi_1,\xi')$ is strictly monotone in $\xi_1$ on the support of a measurable function $\sigma$ on $\mathbb{R}^n$. Then, for all $x_1\in \mathbb{R}$ and with $x'=(x_2,...,x_n)$, we have that
\vspace{0.3cm}
\item[(i)] if $b$ is strictly monotone then
\begin{equation*}\label{CTi}
\n \sigma(D_x)|b'(t)|^{1/2}e^{ib(t)a(\xi)}\varphi(x_1,x')\n_{L^2(\mathbb{R}_t\times\mathbb{R}^{n-1}_{x'})}\leq (2\pi)^{-n}\int_{\mathbb{R}^n}|\widehat{\varphi}(\xi)|^2\frac{|\sigma(\xi)|^2}{|\partial_{\xi_1}a(\xi)|}d\xi.
\end{equation*}
Moreover, if $b$ is strictly monotone and $\lim_{t\rightarrow \infty}b(t)=\infty$, we have
\begin{equation*}\label{CTi2}
\n \sigma(D_x)|b'(t)|^{1/2}e^{ib(t)a(\xi)}\varphi(x_1,x')\n_{L^2(\mathbb{R}_t\times\mathbb{R}^{n-1}_{x'})}=(2\pi)^{-n}\int_{\mathbb{R}^n}|\widehat{\varphi}(\xi)|^2\frac{|\sigma(\xi)|^2}{|\partial_{\xi_1}a(\xi)|}d\xi;
\end{equation*}

\item[(ii)] If $b'$ is such that $\sharp\{t\in\mathbb{R}; b'(t)=0\}=k\geq 1$, then  
\begin{equation*}\label{CTii}
\n \sigma(D_x)|b'(t)|^{1/2}e^{ib(t)a(\xi)}\varphi(x_1,x')\n_{L^2(\mathbb{R}_t\times\mathbb{R}^{n-1}_{x'})}\leq \sqrt{(k+1)}(2\pi)^{-n}\int_{\mathbb{R}^n}|\widehat{\varphi}(\xi)|^2\frac{|\sigma(\xi)|^2}{|\partial_{\xi_1}a(\xi)|}d\xi;
\end{equation*}
\item[(iii)] If $b'$ is such that there exists an increasing sequence of positive times $\{t_k\}_{k\in\mathbb{N}}$ and a decresing sequence of negative times $\{t'_k\}_{k\in\mathbb{N}}$ such that $b'(t_k)=b'(t'_k)=0$, and we set $t_0=t'_0=0$, with $b'$ possibly vanishing at $0$ as well, then, for any function $c\in C(\mathbb{R})$ such that $\sum_{k=0}^\infty \sup_{[t_k,t_{k+1})}|c(t)|$, $\sum_{k=0}^\infty \sup_{(t'_{k+1},t'_{k}]}|c(t)|\leq C<\infty$, we have
\begin{equation*}\label{CTii}
\n \sigma(D_x)|c(t)b'(t)|^{1/2}e^{ib(t)a(\xi)}\varphi(x_1,x')\n_{L^2(\mathbb{R}_t\times\mathbb{R}^{n-1}_{x'})}\leq \sqrt{2C}(2\pi)^{-n}\int_{\mathbb{R}^n}|\widehat{\varphi}(\xi)|^2\frac{|\sigma(\xi)|^2}{|\partial_{\xi_1}a(\xi)|}d\xi.
\end{equation*}

\end{theorem}

\proof
The proof follows  by Lemma \ref{T1} and by application of Theorem 2.1 in \cite{RS} giving the desired inequalities in the case $b(t)=t$.
\endproof

In the sequel we shall make comparisons between the operators  $\LL_{b,a}$ and $\LL_{f,\tilde{a}}$, where, in particular, $f$ satisfies the following assumption (H):
\vspace{0.3cm}
\item[\hypertarget{H}{\bf(H)}] A function $f\in C^1(\mathbb{R})$ is said to satisfy condition (H) if 
\begin{itemize}
\item[-] $f(0)=0$;
\item[-] $f$ is strictly monotone;
\item[-] $\lim_{t\rightarrow \infty}f(t)=\infty$.
\end{itemize}

\begin{corollary}\label{corTC2}
Let $b\in C^1(\mathbb{R)}$ be such that $b(0)=0$ and let $a,\tilde{a}\in C^1(\mathbb{R}^n)$ be real-valued functions such that, for almost all $\xi'=(\xi_2,...,\xi_n)\in\mathbb{R}^n$ we have that $a(\xi_1,\xi'), \tilde{a}(\xi_1,\xi')$ are strictly monotone in $\xi_1$ on the support of a measurable function $\chi$ on $\mathbb{R}^n$. Let $\sigma,\tau\in C^0(\mathbb{R}^n)$ be such that, for some $A>0$
\begin{equation}\label{cond}
\frac{|\sigma(\xi)|}{|\partial_{\xi_1}a(\xi)|^{1/2}}\leq A \frac{|\tau(\xi)|}{|\partial_{\xi_1}\tilde{a}(\xi)|^{1/2}},
\end{equation}
for all $\xi\in\mathsf{supp}\chi$ such that $\partial_{\xi_1}a\neq 0$ and $\partial_{\xi_1}\tilde{a}\neq 0$. Then, for all $f\in C^1(\mathbb{R})$ satisfying condition \hyperlink{H}{(H)}, we have
\vspace{0.3cm}

\item[(i)] if $b$ satisfies (i) or (ii) of Lemma \ref{T1}, then
$$\n|b'(t)|^{1/2}\chi(D_x) \sigma(D_x) e^{i b(t)a(D_x)}\varphi(x_1,x')\n_{L^2(\mathbb{R}_t\times\mathbb{R}^{n-1}_{x'})}$$
$$\leq C A \n  |f'(t)|^{1/2} \chi(D_x) \tau(D_x) e^{i f(t) \tilde{a}(D_x)}\varphi(\tilde{x}_1,x')\n_{L^2(\mathbb{R}_t\times\mathbb{R}^{n-1}_{x'})},$$
where $C=1$ or $C=\sqrt{k+1}$ if $b$ satisfies (i) or (ii) of Lemma \ref{T1}, respectively. Consequently, for any measurable function $\omega$ on $\mathbb{R}$, we have
\begin{equation}
\n |b'(t)|^{1/2}\omega(x_1)\chi(D_x) \sigma(D_x) e^{i b(t) a(D_x)}\varphi(x)\n_{L^2(\mathbb{R}_t\times\mathbb{R}^{n}_{x})}
\end{equation}
$$\leq CA \n  |f'(t)|^{1/2} \omega(x_1)\chi(D_x) \tau(D_x) e^{i f(t) a'(D_x)}\varphi(x)\n_{L^2(\mathbb{R}_t\times\mathbb{R}^{n}_{x})}.$$
Finally, if \eqref{cond} holds with equality and if $b$ satisfies \hyperlink{H}{(H)}, then we get equalities in the previous relations.

\item[(ii)]  if $b$ and $c$ satisfy (iii) of Lemma \ref{T1}, then
$$\n|c(t) b'(t)|^{1/2}\chi(D_x) \sigma(D_x) e^{i b(t)a(D_x)}\varphi(x_1,x')\n_{L^2(\mathbb{R}_t\times\mathbb{R}^{n-1}_{x'})}$$
$$\leq C A \n  |f'(t)|^{1/2} \chi(D_{\tilde{x}_1,x'}) \tau(D_x) e^{i f(t) a'(D_x)}\varphi(\tilde{x}_1,x')\n_{L^2(\mathbb{R}_t\times\mathbb{R}^{n-1}_{x'})},$$
where $C$ depends on the properties of $b$ and $c$. Consequently, for any measurable function $\omega$ on $\mathbb{R}$, we have
\begin{equation}
\n |c(t)b'(t)|^{1/2}\omega(x_1)\chi(D_x) \sigma(D_x) e^{i b(t) a(D_x)}\varphi(x)\n_{L^2(\mathbb{R}_t\times\mathbb{R}^{n}_{x})}
\end{equation}
$$\leq CA \n  |f'(t)|^{1/2} \omega(x_1)\chi(D_x) \tau(D_x) e^{i f(t) a'(D_x)}\varphi(x)\n_{L^2(\mathbb{R}_t\times\mathbb{R}^{n}_{x})}.$$
\vspace{0.05cm}

Conversely, if $\chi\in C^0(\mathbb{R}^n)$, $\omega\neq 0$ on a set of positive measure, $b$ satisfies condition  \hyperlink{H}{(H)} and one of the estimates in (i) is satisfied for all $\varphi$, for some $x_1,\tilde{x}_1\in \mathbb{R}$, then \eqref{cond} holds.
\end{corollary}

\proof
Applying Theorem \ref{TC2} first, then using \eqref{cond} and finally the equality given by point (i) of  Theorem \ref{TC2} (holding when \hyperlink{H}{(H)} is satisfied) the proof follows.
\endproof

From the previous corollary we immediately obtain the following results for $n=1$ and $n=2$.

\begin{corollary}\label{cor1}
Let $b\in C^1(\mathbb{R)}$ be such that $b(0)=0$ and let $a,\tilde{a}\in C^1(\mathbb{R})$ be real-valued functions strictly monotone  on the support of a measurable function $\chi$ on $\mathbb{R}$. Let $\sigma,\tau\in C^0(\mathbb{R})$ be such that there exist $A>0$ such that
\begin{equation}\label{cond2}
\frac{|\sigma(\xi)|}{|\partial_{\xi}a(\xi)|^{1/2}}\leq A \frac{|\tau(\xi)|}{|\partial_{\xi}\tilde{a}(\xi)|^{1/2}},
\end{equation}
for all $\xi\in\mathsf{supp}\chi$ such that $\partial_{\xi} a\neq 0$ and $\partial_{\xi} \tilde{a}\neq 0$. Then,  for all $f\in C^1(\mathbb{R})$ satisfying condition \hyperlink{H}{(H)}, we have
\vspace{0.3cm}

\item[(i)]Let $n=1$. If $b$ satisfies hypothesis (i) or (ii) of Lemma \ref{T1}, then
$$\n|b'(t)|^{1/2}\chi(D_x) \sigma(D_x) e^{i b(t)a(D_x)}\varphi(x)\n_{L^2(\mathbb{R}_t)}$$
$$\leq C A \n  |f'(t)|^{1/2} \chi(D_x) \tau(D_x) e^{i f(t) a'(D_x)}\varphi(\tilde{x})\n_{L^2(\mathbb{R}_t)},$$
where $C=1$ or $C=\sqrt{k+1}$ if $b$ satisfies (i) or (ii) of Lemma \ref{T1}, respectively. Consequently, for any $n\geq 1$ and for any measurable function $\omega$ on $\mathbb{R}^n$, we have
\begin{equation}
\n |b'(t)|^{1/2}\omega(x)\chi(D_j) \sigma(D_j) e^{i b(t) a(D_j)}\varphi(x)\n_{L^2(\mathbb{R}_t\times\mathbb{R}^{n}_{x})}
\end{equation}
$$\leq CA \n  |f'(t)|^{1/2} \omega(x)\chi(D_j) \tau (D_j) e^{i f(t) a'(D_j)}\varphi(x)\n_{L^2(\mathbb{R}_t\times\mathbb{R}^{n}_{x})}.$$
Finally, if \eqref{cond2} holds with equality and if $b$ satisfies \hyperlink{H}{(H)}, then we get equalities in the previous relations.

\item[(ii)] Let $n=1$. If $b$ and $c$ satisfy hypothesis (iii) of Lemma \ref{T1}, then
$$\n|c(t) b'(t)|^{1/2}\chi(D_x) \sigma(D_x) e^{i b(t)a(D_x)}\varphi(x)\n_{L^2(\mathbb{R}_t)}$$
$$\leq C A \n  |f'(t)|^{1/2} \chi(D_{x}) \tau(D_x) e^{i f(t) a'(D_x)}\varphi(\tilde{x})\n_{L^2(\mathbb{R}_t)},$$
where $C$ depends on the properties of $b$ and $c$. Consequently, for any $n\geq 1$ and for any measurable function $\omega$ on $\mathbb{R}^n$, we have
\begin{equation}
\n |c(t)b'(t)|^{1/2}\omega(x)\chi(D_j) \sigma(D_j) e^{i b(t) a(D_j)}\varphi(x)\n_{L^2(\mathbb{R}_t\times\mathbb{R}^{n}_{x})}
\end{equation}
$$\leq CA \n  |f'(t)|^{1/2} \omega(x)\chi(D_j) \ \tau(D_j) e^{i f(t) a'(D_j)}\varphi(x)\n_{L^2(\mathbb{R}_t\times\mathbb{R}^{n}_{x})}.$$
\vspace{0.05cm}

Conversely, if $\chi\in C^0(\mathbb{R})$ and $\omega\neq 0$ on a set of positive measure, $b$ satisfies condition  \hyperlink{H}{(H)}, and one of the estimates in (i) is satisfied for all $\varphi$ and the norms are finite, then inequality \eqref{cond2} holds.
\end{corollary}

\begin{corollary}\label{cor2}
Let $b\in C^1(\mathbb{R)}$ be such that $b(0)=0$ and let $a,\tilde{a}\in C^1(\mathbb{R}^2)$ be real-valued functions strictly monotone  on the support of a measurable function $\chi$ on $\mathbb{R}^2$. Let $\sigma,\tau\in C^0(\mathbb{R}^2)$ be such that, for some $A>0$, we have
\begin{equation}\label{cond3}
\frac{|\sigma(\xi)|}{|\partial_{\xi_1}a(\xi)|^{1/2}}\leq A \frac{|\tau(\xi)|}{|\partial_{\xi_1}\tilde{a}(\xi)|^{1/2}},\quad \xi\in\mathbb{R}^2,
\end{equation}
for all $\xi\in\mathsf{supp}\chi$ such that $\partial_{\xi_1}a\neq 0$ and $\partial_{\xi_1}\tilde{a}\neq 0$. Then,  for all $f\in C^1(\mathbb{R})$ satisfying condition \hyperlink{H}{(H)}, we have
\vspace{0.3cm}

\item[(i)] Let $n=2$. If $b$ satisfies hypothesis (i) or (ii) of Lemma \ref{T1}, then
$$\n|b'(t)|^{1/2}\chi(D_{x_1},D_{x_2}) \sigma(D_{x_1},D_{x_2}) e^{i b(t)a(D_{x_1},D_{x_2})}\varphi(x)\n_{L^2(\mathbb{R}_t\times\mathbb{R}_{x_2})}$$
$$\leq CA \n  |f'(t)|^{1/2} \chi(D_{x_1},D_{x_2}) \tau (D_{x_1},D_{x_2}) e^{i f(t) a'(D_{x_1},D_{x_2})}\varphi(\tilde{x}_1,x_2)\n_{L^2(\mathbb{R}_t\times\mathbb{R}_{x_2})},$$
where $C=1$ or $C=\sqrt{k+1}$ if $b$ satisfies (i) or (ii) of Lemma \ref{T1}, respectively. Consequently, for any $n\geq 1$ and for any measurable function $\omega$ on $\mathbb{R}^{n-1}$, we have
\begin{equation*}
\n |b'(t)|^{1/2}\omega(\check{x}_k)\chi(D_j,D_k) \sigma(D_j,D_k) e^{i b(t) a(D_j,D_k)}\varphi(x)\n_{L^2(\mathbb{R}_t\times\mathbb{R}^{n}_{x})}
\end{equation*}
$$\leq CA \n  |f'(t)|^{1/2} \omega(\check{x}_k)\chi(D_j,D_k) \tau(D_j,D_k) e^{i f(t) a'(D_j,D_k)}\varphi(x)\n_{L^2(\mathbb{R}_t\times\mathbb{R}^{n}_{x})},$$
where $\check{x}_k=(x_1,...,x_{k-1}, x_{k+1},...,x_n)$. Finally, if \eqref{cond3} holds with equality and if $b$ satisfies \hyperlink{H}{(H)}, then we get equalities in the previous relations.
\\
\item[(ii)]  Let $n=2$. If $b$ and $c$ satisfy hypothesis (iii) of Lemma \ref{T1}, then
$$\n|c(t)b'(t)|^{1/2}\chi(D_{x_1},D_{x_2}) \sigma(D_{x_1},D_{x_2}) e^{i b(t)a(D_{x_1},D_{x_2})}\varphi(x)\n_{L^2(\mathbb{R}_t\times\mathbb{R}_{x_2})}$$
$$\leq CA \n  |f'(t)|^{1/2} \chi(D_{x_1},D_{x_2}) \tau (D_{x_1},D_{x_2}) e^{i f(t) a'(D_{x_1},D_{x_2})}\varphi(\tilde{x}_1,x_2)\n_{L^2(\mathbb{R}_t\times\mathbb{R}_{x_2})},$$
where $C$ depends on the properties of $b$ and $c$. Consequently, for any $n\geq 1$ and for any measurable function $\omega\in\mathbb{R}^n$, we have
\begin{equation*}
\n |c(t)b'(t)|^{1/2}\omega(\check{x}_k)\chi(D_j,D_k) \sigma(D_j,D_k) e^{i b(t) a(D_j,D_k)}\varphi(x)\n_{L^2(\mathbb{R}_t\times\mathbb{R}^{n}_{x})}
\end{equation*}
$$\leq CA \n  |f'(t)|^{1/2} \omega(\check{x}_k)\chi(D_j,D_k) \tau(D_j,D_k) e^{i f(t) a'(D_j,D_k)}\varphi(x)\n_{L^2(\mathbb{R}_t\times\mathbb{R}^{n}_{x})},$$
where $\check{x}_k=(x_1,...,x_{k-1}, x_{k+1},...,x_n)$.\vspace{0.3cm}

Conversely, if $\chi\in C^0(\mathbb{R}^2)$, $\omega\neq 0$ on a set of $\mathbb{R}^{n-1}$ with positive measure, $b$ satisfies condition \hyperlink{H}{(H)} and one of the estimates in (i) is satisfied for all $\varphi$, then inequality \eqref{cond3} holds.
\end{corollary}

\begin{remark}
Notice that, by virtue of the previous results, one can also compare time-degenerate operators with different degeneracies (i.e. $\LL_{b,a}$ and $\LL_{f,\tilde{a}}$), provided that $b$ and $f$ satisfy suitable conditions. This will therefore allow comparisons between the standard Schr\"{o}dinger operator $\LL_{t,\Delta}$ and time-degenerate  Schr\"{o}dinger operators of the form $\LL_{b,\Delta}$.
\end{remark}
\vspace{0.3cm}

\noindent {\bf The radially symmetric case.} 
In the radially symmetric case comparison principles similar to the previous ones can be obtained. Below we derive the suitable generalization to our case of the comparison principle holding in the nondegenerate radially symmetric case proved in \cite{RS} .

\begin{theorem}\label{rad}
Let $a, \tilde{a}\in C^1(\mathbb{R}_+)$ be real-valued strictly monotone on the support of a measurable function $\chi$ on $\mathbb{R}_+$, and let  $b,f\in C^1(\mathbb{R})$ vanish at 0, with $f$ satisfying condition \hyperlink{H}{(H)}. Let also $\sigma,\tau\in C^0(\mathbb{R}_+)$ be such that

\begin{equation}\label{cond4}
\frac{|\sigma(\rho)|}{|\frac{d}{d\rho}a(\rho)|^{1/2}}\leq A \frac{|\tau(\xi)|}{|\frac{d}{d\rho}\tilde{a}(\rho)|^{1/2}},
\end{equation}
for all $\rho \in \mathsf{supp}\chi$ where $\frac{d}{d\rho}a(\rho)\neq 0$ and $\frac{d}{d\rho}\tilde{a}(\rho)\neq 0$. Then we have
\vspace{0.3cm}

\item[(i)] If $b$ satisfies hypothesis (i) or (ii) of Lemma \ref{T1} then
\begin{equation}\label{sym1}
\n \chi(|D_x|)|b'(t)|^{1/2}\sigma(|D_x|))e^{ib(t)a(|D_x|)}\varphi(x)\n_{L^2(\mathbb{R}_t\times\mathbb{R}^n_x)} 
\end{equation}
$$\leq CA \n \chi(|D_x|)|f'(t)|^{1/2} \tau(|D_x|))e^{if(t) \tilde{a}(|D_x|)}\varphi(x)\n_{L^2(\mathbb{R}_t\times\mathbb{R}^n_x)},
$$
where $C=1$ or $C=\sqrt{(k+1)}$ if $b$ satisfies (i) or (ii) of Lemma \ref{T1} respectively.

Moreover, for any measurable function $\omega$ on $\mathbb{R}^n$, we have
\begin{equation}\label{sym1.2}
\n \omega(x)\chi(|D_x|)|b'(t)|^{1/2}\sigma(|D_x|))e^{ib(t)a(|D_x|)}\varphi(x)\n_{L^2(\mathbb{R}_t\times\mathbb{R}^n_x)} 
\end{equation}
$$\leq CA \n \omega(x) \chi(|D_x|)|f'(t)|^{1/2} \tau(|D_x|))e^{if(t) \tilde{a}(|D_x|)}\varphi(x)\n_{L^2(\mathbb{R}_t\times\mathbb{R}^n_x)},
$$
where $C$ is the constant in \eqref{sym1}. Finally, if \eqref{cond4} holds with equality and if $b$ satisfies \hyperlink{H}{(H)}, then we get equalities in the previous relations.
\\
\item[(ii)]  If $b$ and $c$ satisfy hypothesis (iii) of Lemma \ref{T1} then 
\begin{equation*}\label{sym2}
\n \chi(|D_x|)|c(t)b'(t)|^{1/2}\sigma(|D_x|))e^{ib(t)a(|D_x|)}\varphi(x)\n_{L^2(\mathbb{R}_t\times\mathbb{R}^n_x)} 
\end{equation*}
$$\leq C' A \n \chi(|D_x|)|f'(t)|^{1/2} \tau(|D_x|))e^{if(t) \tilde{a}(|D_x|)}\varphi(x)\n_{L^2(\mathbb{R}_t\times\mathbb{R}^n_x)},
$$
where $C'$ depends on $b$ and $c$. 

Moreover, for any measurable function $\omega$ on $\mathbb{R}^n$, we have
\begin{equation*}
\n \omega(x)\chi(|D_x|)|c(t)b'(t)|^{1/2}\sigma(|D_x|))e^{ib(t)a(|D_x|)}\varphi(x)\n_{L^2(\mathbb{R}_t\times\mathbb{R}^n_x)} 
\end{equation*}
$$\leq C'A \n \omega(x) \chi(|D_x|)|f'(t)|^{1/2} \tau(|D_x|))e^{if(t) \tilde{a}(|D_x|)}\varphi(x)\n_{L^2(\mathbb{R}_t\times\mathbb{R}^n_x)},
$$
where $C'$ is the same constant as above.
\vspace{0.3cm}

Conversely, if $\chi\in C^0(\mathbb{R}_+)$, $\omega\neq 0$ on a set of $\mathbb{R}^n$ with positive measure, $b$ satisfies condition \hyperlink{H}{(H)} and \eqref{sym1} (for some $x,\tilde{x}$) or  \eqref{sym1.2} is satisfied for any $\varphi$ and the norms are finite, then inequality \eqref{cond4} holds.
\end{theorem}
\proof
The proof follows by  Lemma \ref{T1} and Theorem 2.5 in \cite{RS}. The converse part holds only if $b$ satisfies \hyperlink{H}{(H)}, since, in this case,
$$\n \chi(|D_x|)|b'(t)|^{1/2}\sigma(|D_x|))e^{ib(t)a(|D_x|)}\varphi(x)\n_{L^2(\mathbb{R}_t\times\mathbb{R}^n_x)}$$
$$=\n \chi(|D_x|)\sigma(|D_x|))e^{ita(|D_x|)}\varphi(x)\n_{L^2(\mathbb{R}_t\times\mathbb{R}^n_x)},$$
and the result follows again from Theorem 2.5 in \cite{RS}.
\endproof

\section{Application of comparison principles to smoothing estimates}\label{ACP}
In this section we will will show how to use comparison principles to obtain smoothing estimates for $\LL_{b,a}$ provided that it is comparable to an operator $\LL_{f,\tilde{a}}$, with   $\LL_{f,\tilde{a}}$ satisfying smoothing estimates.

In addition, if the operators $\LL_{b,a}$ and $\LL_{f,\tilde{a}}$ are equivalent in some sense, that is, if conditions of the form \eqref{cond} hold in both directions, then we have an equivalence between the compared quantities, and one has that one operator satisfies a smoothing estimates if and only if the other one satisfies the suitable corresponding smoothing estimate.  In this case, as we shall see below, in order to obtain equivalent norms we will need to restrict ourselves to the case when $b,f\in C^1({\mathbb{R}})$ and are such that they satisfy condition \hyperlink{H}{(H)}, since, otherwise, the equivalence given by \eqref{cond} would not be enough to get the equivalence of the associated weighted norms.

Here we shall specify the analysis to the case $a(D)=|D_x|^m$ and $\tilde{a}(D)=|D_x|^l$, with $|D_x|$ the Fourier multiplier with symbol $|\xi|$, to show first that weighted smoothing estimates can be derived by comparison, and, afterwards, that smoothing estimates are all equivalent, or, more precisely, that weighted smoothing estimates for different operators can be derived by equivalence.
\vspace{0.3cm}

We start by showing, by comparison, some smoothing estimates in the cases $n=1$ and $n=2$ ($n$ is the space-dimension). When $n=1$ we consider $\LL_{b,|D_x|^m}$ and $\LL_{f, |D_x|^l}$, with
$l,m>0$, and we set $a(\xi)=|\xi|^m$ and $\tilde{a}(\xi)=|\xi|^l$. We then consider $\sigma(\xi)=|\xi|^{\frac{m-1}{2}}$ and $\tau(\xi)=|\xi|^{\frac{l-1}{2}}$ and have that \eqref{cond} is satisfied, since
$$\frac{|\sigma(\xi)|}{|\frac{d}{d\xi}a(\xi)|^{1/2}}=\sqrt{\frac l m} \frac{|\tau(\xi)|}{|\frac{d}{d\xi}\tilde{a}(\xi)|^{1/2}}.$$ 
When $n=2$ we can make similar choices for the quantities above, and have, again, that \eqref{cond} is satisfied. 

It is important to observe that we are allowed to consider $f=t$ and $l=1$. This choice, in combination with the estimate $\n e^{it|D_x|}\varphi\n_{L^2_t}=\n\varphi\n_{L^2_x}$, gives the following theorem about the homogeneous smoothing effect for time-degenerate Schr\"{o}dinger-type operators.

\begin{theorem}\label{Thm2}
Let $m>0$, and let $b\in C^1(\mathbb{R})$ be such that it vanishes at 0. Then, for all $x\in\mathbb{R}^n$ and for all $\varphi\in L^2(\mathbb{R}^n)$, we have
\vspace{0.3cm}

\item[(i)] If $n=1$ and $b$ satisfies (i) or (ii) of Lemma \ref{T1}, we have
$$\sup_{x\in\mathbb{R}}\n |b'(t)|^{1/2}|D_{x}|^{\frac{m-1}{2}} e^{ib(t)|D_{x}|^{m}}\varphi\n_{L^2(\mathbb{R}_{t})}
\leq C \n \varphi\n_{L^2(\mathbb{R}_x)},$$
and, if $n=2$,
$$\sup_{x_1\in\mathbb{R}}\n |b'(t)|^{1/2} |D_{x_2}|^{\frac{m-1}{2}} e^{ib(t)D_{x_1}|D_{x_2}|^{m-1}}\varphi\n_{L^2(\mathbb{R}_{x_2}\times \mathbb{R}_t)}
\leq C \n \varphi\n_{L^2(\mathbb{R}^2_x)}.$$
\item[(ii)] If $n=1$ and $b$ and $c$ satisfy (iii) of Lemma \ref{T1}, we have
$$ \sup_{x\in\mathbb{R}}\n |c(t)b'(t)|^{1/2}|D_{x}|^{\frac{m-1}{2}} e^{ib(t)|D_{x}|^{m}}\varphi\n_{L^2(\mathbb{R}_{t})}
\leq C \n \varphi\n_{L^2(\mathbb{R}_x)},$$
and, if $n=2$,
$$\sup_{x_1\in\mathbb{R}}\n |c(t)b'(t)|^{1/2}|D_{x_2}|^{\frac{m-1}{2}} e^{ib(t)D_{x_1}|D_{x_2}|^{m-1}}\varphi\n_{L^2(\mathbb{R}_{x_2}\times \mathbb{R}_t)}
\leq C \n \varphi\n_{L^2(\mathbb{R}^2_x)}.$$
\end{theorem}
\proof
The proof follows by application of comparison principles. 
As regards point (i) when $n=1$, denoting by $\chi_\pm:=\chi_{\mathbb{R}_\pm}$, we have that 
$$\sup_x\n |b'(t)|^{1/2}|D_x|^{\frac{m-1}{2}}e^{i b(t) |D_x|^m}\varphi\n_{L^2_t(\mathbb{R})}$$
$$=\sup_x\n (\chi_+(D_x)+\chi_-(D_x))|b'(t)|^{1/2}|D_x|^{\frac{m-1}{2}}e^{i b(t) |D_x|^m}\varphi\n_{L^2_t(\mathbb{R})}$$
$$\leq \sup_x\n |b'(t)|^{1/2}\chi_+(D_x)|D_x|^{\frac{m-1}{2}}e^{i b(t) |D_x|^m}\varphi\n_{L^2_t(\mathbb{R})}$$
$$+\sup_x\n |b'(t)|^{1/2} \chi_+(D_x)|D_x|^{\frac{m-1}{2}}e^{i b(t) |D_x|^m}u_0\n_{L^2_t(\mathbb{R})}.$$
Now, by point (i) of Corollary \ref{cor1} with  $f(t)=t$, $\tilde{a}(D_x)=|D_x|$ and $\tau(\xi)=1$, together with the estimate $\n e^{it|D_x|}\varphi\n_{L^2_t}=\n\varphi\n_{L^2_x}$, we get
$$ \sup_x\n |b'(t)|^{1/2}\sigma_\pm(D_x)|D_x|^{\frac{m-1}{2}}e^{i b(t) |D_x|^m}\varphi\n_{L^2_t(\mathbb{R})}
 \leq C \n \varphi\n_{L^2_x(\mathbb{R}},$$
 which proves point (i) when $n=1$. 
 
 As regards point (i) when $n=2$, it follows by application by point (i) of Corollary \ref{cor2} with $f=t$, $\tilde{a}(\xi)=\xi_1$ and $\tau(\xi)=1$. In fact, denoting by $\chi_\pm=\chi_{\mathbb{R}_{{x_2}_\pm}}$, we have
$$ \n (\chi_+(D_{x_2})+\chi_-(D_{x_2}))|b'(t)|^{1/2} |D_{x_2}|^{\frac{m-1}{2}} e^{ib(t)D_{x_1}|D_{x_2}|^{m-1}}\varphi\n_{L^2(\mathbb{R}_{x_2}\times \mathbb{R}_t)}\leq C \n e^{t|D_{x_1}|}\varphi\n_{L^2(\mathbb{R}_{x_2}\times \mathbb{R}_t)}$$
$$=C\n \varphi\n_{L^2_x},$$
 which gives the result.
  
 The proof of point (ii) reads exactly the same, the only difference is that we use point (ii) of Corollary \ref{cor1} and point (ii) of Corollary \ref{cor2} when $n=1$ and $n=2$ respectively.
 \endproof
Observe that Theorem \ref{Thm2} applied to the case $b(t)=t$ gives the standard results proved by Linares and Ponce in \cite{LP}.

\begin{remark}
Note that when ${\bf{n=1}}$, by (i) of Corollary \ref{cor1}, if $f$ and $b$ satisfy \hyperlink{H}{(H)}, then, for all $\varphi$ such that either $\mathsf{supp}\, \widehat{\varphi} \subset [0,+\infty)$ or $\mathsf{supp}\, \widehat{\varphi} \subset (-\infty,0]$ (i.e. $a(\xi)$ and $\tilde{a}(\xi)$ are strictly monotone on $\mathsf{supp} \, \widehat{\varphi}$), we have
\begin{equation*} \label{eqn1}
\n |D_x|^{\frac{m-1}{2}} |b'(t)|^{1/2}e^{ib(t)|D_x|^m}\varphi\n_{L^2_t(\mathbb{R})}=\sqrt{\frac l m}\n |D_x|^{\frac{l-1}{2}} |f'(t)|^{1/2}e^{if(t)|D_x|^l}\varphi\n_{L^2_t(\mathbb{R})},\quad \forall x\in\mathbb{R},
\end{equation*}
where, of course, the case $f=b$ is possible.

In the case ${\bf{n=2}}$, under the same hypotheses on the functions $b$ and $f$ (possibly $b=f$) as before, we can apply point (i) of Corollary \ref{cor2} and get, for all $\varphi$ such that either $\pi_2(\mathsf{supp}\, \widehat{\varphi}) \subset [0,+\infty)$ or $\pi_2(\mathsf{supp}\, \widehat{\varphi}) \subset (-\infty,0]$ (where $\pi_2$ is the canonical projection on $\mathbb{R}_{\xi_2})$,
\begin{equation*} \label{eqn2}
\n |D_{x_2}|^{\frac{m-1}{2}} |b'(t)|^{1/2}e^{ib(t)D_{x_1}|D_{x_2}|^{m-1}}\varphi\n_{L^2(\mathbb{R}_{x_2}\times \mathbb{R}_t)}$$
$$=\n |D_{x_2}|^{\frac{l-1}{2}} |f'(t)|^{1/2}e^{if(t)D_{x_1}|D_{x_2}|^{\frac{l-1}{2}}}\varphi\n_{L^2(\mathbb{R}_{x_2}\times \mathbb{R}_t)},
\end{equation*}
for all $x\in\mathbb{R}^2$ and for all $\varphi$ (for which the previous norms are well defined), where, again, the case $f=b$ is possible.

These properties show that, if the two operators are {\it equivalent} in some sense, that is if the weighted norms are equivalent, then one operator satisfies a smoothing estimate if and only if the other one satisfies the suitable corresponding smoothing estimate. 
An equivalence relation between the norms is established by the properties of the functions $b$ and $f$ and by relation \eqref{cond}.
In particular, due to the previous identities, we can say that the weighted norms associated with $\LL_{b,|D_x|^m}$ and  $\LL_{f,|D_x|^l}$ are {\it{equivalent} } if \eqref{cond} holds with equality and if both $b$ and $f$ satisfy condition \hyperlink{H}{(H)}.
\end{remark}

Point (i) of Theorem \ref{Thm2} applied to the case $n=1$, $m=2$, gives the natural extension of the standard homogeneous smoothing effect holding for Schr\"{o}dinger operators. 
By using comparison principles we can also obtain the following multidimensonal  version.

\begin{theorem}\label{GSE}
Let $n\geq 1$ and $b\in C^1(\mathbb{R})$ be such that it vanishes at 0. Then, for all $x\in\mathbb{R}^n$ and for all $\varphi\in L^2(\mathbb{R}^n)$, 
\vspace{0.3cm}

\item[(i)] If $b$ satisfies (i) or (ii) of Lemma \ref{T1}, we have
$$\sup_{x_j}\n |b'(t)|^{1/2}|D_{x_j}|^{\frac{1}{2}} e^{ib(t)\Delta_x}\varphi\n_{L^2(\mathbb{R}_{x'}\times \mathbb{R}_t)}
\leq C \n \varphi\n_{L^2(\mathbb{R}_x)},$$
where $x'=(x_1,...,x_{j-1},x_{j+1},...,x_n)$.

\item[(ii)] If $b$ and $c$ satisfy (iii) of Lemma \ref{T1}, we have
$$\sup_{x_j}\n |c(t)b'(t)|^{1/2}|D_{x_j}|^{\frac{1}{2}} e^{ib(t)\Delta_x}\varphi\n_{L^2(\mathbb{R}_{x'}\times \mathbb{R}_t)}
\leq C \n \varphi\n_{L^2(\mathbb{R}_x)},$$
where $x'=(x_1,...,x_{j-1},x_{j+1},...,x_n).$
\end{theorem}
\proof
The proof follows by application of the procedure used in Theorem \ref{Thm2}.
\endproof

We shall now derive other smoothing type estimates by using results obtained by Sugimoto in \cite{S} and by Walther in \cite{W} combined with the comparison principles obtained above. In particular it was proved by Sugimoto in \cite{S} that, given $n\geq 2$ and $1-n/2<\beta<1/2$, for all $\varphi\in L^2(\mathbb{R}^n)$ we have 
\begin{equation}
\label{sug}
\n |x|^{\beta-1}|D_x|^\beta e^{it|D_x|^2}\varphi\n_{L^2(\mathbb{R}_t\times\mathbb{R}^n_x)}\leq C\n\varphi\n_{L^2(\mathbb{R}^n_x)}.
\end{equation}
Notice that the previous result is a consequence of the result by Kato and Yajima in \cite{Kato-Yajima} showing that, for all $\varphi\in L^2(\mathbb{R}^n)$
\begin{equation}
\label{ky}
\n |x|^{\beta-1}|D_x|^\beta e^{it|D_x|^2}\varphi\n_{L^2(\mathbb{R}_t\times\mathbb{R}^n_x)}\leq C\n\varphi\n_{L^2(\mathbb{R}^n_x)},\quad 1/2-\varepsilon\leq \beta<1/2,
\end{equation}
with $0<\varepsilon<1/2$. To see this it is enough to observe that, by Theorem $B^*$ in \cite{SW},
$$\n |x|^{\beta-1}|D_x|^\beta e^{it|D_x|^2}\varphi\n_{L^2(\mathbb{R}_t\times\mathbb{R}^n_x)}\leq
\n |x|^{(1/2-\varepsilon)-1}|D_x|^{1/2-\varepsilon}e^{it|D_x|^2}\varphi\n_{L^2(\mathbb{R}_t\times\mathbb{R}^n_x)}$$
$$\underset{\eqref{ky}}{\leq} C\n\varphi\n_{L^2(\mathbb{R}^n_x)}.$$

As mentioned before, we will also make use of a result by Walther (Theorem 4.1 in \cite{W}) in which he shows that, for $n,m>1$, and for all $\varphi\in L^2(\mathbb{R}^n)$ such that $\mathsf{supp}\widehat{\varphi}\subset \{\xi\in\mathbb{R}^n:\,|\xi|\leq 1\}$, the following estimate holds
\begin{equation}
\label{w}
\n \la x\ra^{-m/2}e^{it|D_x|^{m}}\varphi\n_{L^2(\mathbb{R}_t\times\mathbb{R}^n_x)}\leq C\n\varphi\n_{L^2(\mathbb{R}^n_x)}.
\end{equation}
We can then prove the following theorem.

\begin{theorem}
Let $n\geq 2$ and let $b,f\in C^1(\mathbb{R})$ be such that they satisfy condition \hyperlink{H}{(H)}. Then, for all $\varphi\in L^2(\mathbb{R}^n)$, we have
\begin{equation}
\label{sugb}
\n |x|^{\beta-1}|b'(t)|^{1/2}|D_x|^\beta e^{ib(t)|D_x|^2
}\varphi\n_{L^2(\mathbb{R}_t\times\mathbb{R}^n_x)}
\leq C\n\varphi\n_{L^2(\mathbb{R}^n_x)},
\end{equation}
for $1-n/2<\beta<1/2$, and
\begin{equation}
\label{sugf}
\n |x|^{\alpha-m/2}|f'(t)|^{1/2}|D_x|^{\alpha}e^{if(t)|D_x|^{m}}\varphi\n_{L^2(\mathbb{R}_t\times\mathbb{R}^n_x)}\leq C\n\varphi\n_{L^2(\mathbb{R}^n_x)},
\end{equation}
for $(m-n)/2<\alpha<(m-1)/2$.

Moreover \eqref{sug}, \eqref{sugf}, \eqref{sugb}, \eqref{w} and \eqref{ky} (with  $\varepsilon>0$ sufficiently small) are equivalent.
Finally, for $m>0$ and for any $\alpha,\beta>0$, we have, for all $\varphi\in L^2(\mathbb{R}^n)$ such that $\mathsf{supp}\widehat{\varphi}\subset \{\xi\in\mathbb{R}^n:\,|\xi|\leq 1\}$,
\begin{equation}
\label{1}
\n |x|^{\beta-1}|b'(t)|^{1/2}|D_x|^\beta e^{ib(t)|D_x|^2}\varphi\n_{L^2(\mathbb{R}_t\times\mathbb{R}^n_x)}=
\end{equation}
$$\sqrt{\frac m 2}\n |x|^{\beta-1}|f'(t)|^{1/2}|D_x|^{m/2+\beta-1}e^{if(t)|D_x|^{m}}\varphi\n_{L^2(\mathbb{R}_t\times\mathbb{R}^n_x)},$$
and, for all $\varphi\in L^2(\mathbb{R}^n)$,
\begin{equation}
\label{2}
\n \la x\ra^{\alpha-m/2}|f'(t)|^{1/2}|D_x|^{\alpha}e^{if(t)|D_x|^{m}}\varphi\n_{L^2(\mathbb{R}_t\times\mathbb{R}^n_x)}
\end{equation}
$$\leq \n |x|^{\alpha-m/2}|f'(t)|^{1/2}|D_x|^{\alpha}e^{if(t)|D_x|^{m}}\varphi\n_{L^2(\mathbb{R}_t\times\mathbb{R}^n_x)}$$
$$\leq \sup_{\lambda>0}\n \la x\ra ^{\alpha-m/2}|f'(t)|^{1/2}|D_x|^{\alpha}e^{if(t)|D_x|^{m}}\varphi_\lambda\n_{L^2(\mathbb{R}_t\times\mathbb{R}^n_x)},$$
where $\varphi_\lambda(x)=\lambda^{n/2}\varphi(\lambda x)$ and $\alpha\leq m/2$ in the last estimate. The previous norms are finite for suitable $\alpha$ and $\beta$ (that is, as in \eqref{sugb},\eqref{sugf}) and the operators $|x|^{\alpha-m/2}|f'(t)|^{1/2}|D_x|^{\alpha}e^{if(t)|D_x|^{m}}$ and $\la x\ra^{\alpha-m/2}|f'(t)|^{1/2}|D_x|^{\alpha}e^{if(t)|D_x|^{m}}$ have the same norms as mapping from $L^2(\mathbb{R}^n_x)$ to $L^2(\mathbb{R}_t\times\mathbb{R}^n_x)$.
\end{theorem}

\proof
Recall that, by Lemma \ref{T1}, if $b\in C^1(\mathbb{R})$ is such that it satisfies \hyperlink{H}{(H)}, then
\begin{equation}\label{inv}
\n |x|^{\beta-1}|b'(t)|^{1/2}|D_x|^\beta e^{ib(t)|D_x|^2}\varphi\n_{L^2(\mathbb{R}_t\times\mathbb{R}^n_x)}=
\n |x|^{\beta-1}|D_x|^\beta e^{it|D_x|^2}\varphi\n_{L^2(\mathbb{R}_t\times\mathbb{R}^n_x)},\end{equation}
therefore, by Theorem 1.1 in \cite{S} (i.e. \eqref{sug}), we also have, for $n\geq 2$, $1-n/2<\beta<1/2$ and $\varphi\in L^2(\mathbb{R}^n)$, 
\begin{equation*}
\n |x|^{\beta-1}|b'(t)|^{1/2}|D_x|^\beta e^{ib(t)|D_x|^2}\varphi\n_{L^2(\mathbb{R}_t\times\mathbb{R}^n_x)}
\leq C\n\varphi\n_{L^2(\mathbb{R}^n_x)},
\end{equation*}
which gives \eqref{sugb} and the equivalence of \eqref{sug} and \eqref{sugb}.

In particular, by point (i) in Theorem \ref{rad}, with $b$ satisfying \hyperlink{H}{(H)}, $m>0$, and $\varphi\in L^2(\mathbb{R}^n)$ such that $\mathsf{supp}\widehat{\varphi}\subset \{\xi\in\mathbb{R}^n:\,|\xi|\leq 1\}$, we have
\begin{equation*}
\n |x|^{\beta-1}|b'(t)|^{1/2}|D_x|^\beta e^{ib(t)|D_x|^2}\varphi\n_{L^2(\mathbb{R}_t\times\mathbb{R}^n_x)}\end{equation*}
$$=\sqrt{\frac m 2}\n |x|^{\beta-1}|f'(t)|^{1/2}|D_x|^{m/2+\beta-1}e^{if(t)|D_x|^{m}}\varphi\n_{L^2(\mathbb{R}_t\times\mathbb{R}^n_x)}$$
$$\underset{(\alpha=m/2+\beta-1)}{=} \sqrt{\frac m 2}\n |x|^{\alpha-m/2}|f'(t)|^{1/2}|D_x|^{\alpha}e^{if(t)|D_x|^{m}}\varphi\n_{L^2(\mathbb{R}_t\times\mathbb{R}^n_x)}$$
(where $f=b$ is allowed) which  gives \eqref{1}. Similarly, by Corollary \ref{cor1}, for all $\varphi\in L^2(\mathbb{R}^n)$ (whose Fourier transform does not necessarily have compact support), we have
 \begin{equation*}
\n |x|^{\beta-1}|b'(t)|^{1/2}|D_x|^\beta e^{ib(t)|D_x|^2}\varphi\n_{L^2(\mathbb{R}_t\times\mathbb{R}^n_x)}\end{equation*}
$$=\left(\n \chi_{\mathbb{R}_+}(D_x) |x|^{\beta-1}|b'(t)|^{1/2}|D_x|^\beta e^{ib(t)|D_x|^2}\varphi\n^2_{L^2(\mathbb{R}_t\times\mathbb{R}^n_x)}\right.$$
$$+\left. \n \chi_{\mathbb{R}_-}(D_x) |x|^{\beta-1}|b'(t)|^{1/2}|D_x|^\beta e^{ib(t)|D_x|^2}\varphi\n^2_{L^2(\mathbb{R}^2_t\times\mathbb{R}^n_x)} \right)^{1/2}$$
$$\leq C \sqrt{\frac m 2}\n |x|^{\alpha-m/2}|f'(t)|^{1/2}|D_x|^{\alpha}e^{if(t)|D_x|^{m}}\varphi\n_{L^2(\mathbb{R}_t\times\mathbb{R}^n_x)}$$
which implies, by \eqref{sugb}, 
\begin{equation*}
\n |x|^{\alpha-m/2}|f'(t)|^{1/2}|D_x|^{\alpha}e^{if(t)|D_x|^{m}}\varphi\n_{L^2(\mathbb{R}_t\times\mathbb{R}^n_x)}\leq C\n\varphi\n_{L^2(\mathbb{R}^n_x)},
\end{equation*}
for $(m-n)/2<\alpha<(m-1)/2$, which gives \eqref{sugf} and the equivalence of \eqref{sugf} and \eqref{sugb}, and, as a consequence, that of \eqref{sugf} and \eqref{sug}.

The equivalence of \eqref{sugb} and \eqref{sugf} with \eqref{ky} is derived from the equivalence of  \eqref{ky} and \eqref{sug}. The latter equivalence is immediate, since, as previously observed, \eqref{ky} implies \eqref{sug} while the opposite is trivial. 

As regards the equivalence with \eqref{w}, on one hand we have that  \eqref{w} is implied by \eqref{sugf} (with $\alpha=0$) together with property \eqref{inv} and the trivial inequality $\la x \ra^{-m/2}\leq | x |^{-m/2}$.
On the other hand, \eqref{w} implies \eqref{sugf}. In fact, given $\chi\in C_0^\infty([0,1))$ such that $\chi(\rho)=1$ for $\rho\leq 1/2$, we have, by comparison,
$$\n \la x\ra^{\alpha-m/2}\chi(|D_x|)|f'(t)|^{1/2}|D_x|^{\alpha}e^{if(t)|D_x|^{m}}\varphi\n_{L^2(\mathbb{R}_t\times\mathbb{R}^n_x)}
\leq$$
$$ \sqrt{\frac \mu m} \n  \la x\ra^{\alpha-m/2}\chi(|D_x|)|f'(t)|^{1/2}|D_x|^{\alpha+(\mu-m)/2}e^{if(t)|D_x|^{\mu}}\varphi\n_{L^2(\mathbb{R}_t\times\mathbb{R}^n_x)}$$
$$\underset{(m-2\alpha=\mu)}{=}\sqrt{\frac{\mu}{m}}\n \la x\ra^{-\mu/2}\chi(|D_x|)|f'(t)|^{1/2}|D_x|^{\alpha}e^{if(t)|D_x|^{\mu}}\varphi\n_{L^2(\mathbb{R}_t\times\mathbb{R}^n_x)}$$
$$=\sqrt{\frac{\mu}{m}} \n \la x\ra^{-\mu/2}\chi(|D_x|)|D_x|^{\alpha}e^{it|D_x|^{\mu}}\varphi\n_{L^2(\mathbb{R}_t\times\mathbb{R}^n_x)},$$
which, by \eqref{w}, gives
$$\n \la x\ra^{\alpha-m/2}\chi(|D_x|)|f'(t)|^{1/2}|D_x|^{\alpha}e^{if(t)|D_x|^{m}}\varphi\n_{L^2(\mathbb{R}_t\times\mathbb{R}^n_x)}
\leq C\n \varphi\n_{L^2(\mathbb{R}^n)}.$$
Finally from the latter we get \eqref{sugf}, since
$$\n |x|^{\alpha-m/2}|f'(t)|^{1/2}|D_x|^{\alpha}e^{if(t)|D_x|^{m}}\varphi\n_{L^2(\mathbb{R}_t\times\mathbb{R}^n_x)}$$
$$=\lim_{\lambda \rightarrow 0^+}\n \lambda^{\alpha-m/2}  \la x/\lambda \ra^{\alpha-m/2} \chi(\lambda|D_x|)|f'(t)|^{1/2}|D_x|^{\alpha}e^{if(t)|D_x|^{m}}\varphi\n_{L^2(\mathbb{R}_t\times\mathbb{R}^n_x)}$$
$$\leq \sup_{\lambda>0}\n \la x \ra^{\alpha-m/2} \chi(|D_x|)|f'(t)|^{1/2}|D_x|^{\alpha}e^{if(t)|D_x|^{m}}\varphi_\lambda\n_{L^2(\mathbb{R}_t\times\mathbb{R}^n_x)}\leq C \sup_{\lambda>0} \n \varphi_\lambda\n_{L^2(\mathbb{R}^n)},$$
with $\varphi_\lambda(x)=\lambda^{n/2}\varphi(\lambda x)$, which gives, in particular, inequality \eqref{sugf} (since $\n \varphi_\lambda\n_{L^2(\mathbb{R}^n)}= \n \varphi\n_{L^2(\mathbb{R}^n)}$). Finally, from the previous inequality, we also obtain \eqref{2}, that is,
$$\n |x|^{\alpha-m/2}|f'(t)|^{1/2}|D_x|^{\alpha}e^{if(t)|D_x|^{m}}\varphi\n_{L^2(\mathbb{R}_t\times\mathbb{R}^n_x)}$$
$$\leq \sup_{\lambda>0}\n \la x \ra^{\alpha-m/2} |f'(t)|^{1/2}|D_x|^{\alpha}e^{if(t)|D_x|^{m}}\varphi_\lambda\n_{L^2(\mathbb{R}_t\times\mathbb{R}^n_x)},$$
which concludes the proof.
\endproof
Notice that all these results coincide with the classical ones when $f(t)=b(t)=t$ (see \cite{KPV2,LP,RS,Kato-Yajima,W}).
From the previous result we derive the following corollaries.

\begin{corollary}
Let $n\geq 1, m>0$, $s>1/2$ and $b\in C^1(\mathbb{R})$ be such that it satisfies \hyperlink{H}{(H)}. Then, for all $\varphi\in L^2(\mathbb{R}^n)$,
\begin{equation*}
\n \la x_n\ra^{-s}|b'(t)|^{1/2}|D_n|^{(m-1)/2}e^{ib(t)|D_n|^{m}}\varphi\n_{L^2(\mathbb{R}_t\times\mathbb{R}^n_x)}\leq \n \varphi\n_{L^2(\mathbb{R}^n)}.
\end{equation*}

Let $n\geq 2, m>0$, $s>1/2$ and $b\in C^1(\mathbb{R})$ be such that it satisfies \hyperlink{H}{(H)}. Then, for all $\varphi\in L^2(\mathbb{R}^n)$
\begin{equation*}
\n \la x_1\ra^{-s}|b'(t)|^{1/2}|D_n|^{(m-1)/2}e^{ib(t)D_1|D_n|^{m-1}}\varphi\n_{L^2(\mathbb{R}_t\times\mathbb{R}^n_x)}\leq \n \varphi\n_{L^2(\mathbb{R}^n)}.
\end{equation*}
\end{corollary}

\begin{corollary}
Let $m>0$, $(m-n+1)/2<\alpha<(m-1)/2$ and $b\in C^1(\mathbb{R})$ be such that it satisfies \hyperlink{H}{(H)}. Then, for all $\varphi\in L^2(\mathbb{R}^n)$,
\begin{equation*}
\n  |x|^{\alpha-m/2}|b'(t)|^{1/2}|D_n|^{(m-1)/2}e^{ib(t)(|D_1|^m-|D'|^{m})}\varphi\n_{L^2(\mathbb{R}_t\times\mathbb{R}^n_x)}\leq \n \varphi\n_{L^2(\mathbb{R}^n)},
\end{equation*}
where $D'=(D_2,...,D_n)$.
\end{corollary}

When $b$ does not satisfy condition \hyperlink{H}{(H)} we can get, by comparison, the same smoothing estimates. However, it is important to stress that we do not have equivalences in this case (i.e. when $b$ does not satisfy condition \hyperlink{H}{(H)}).
\begin{corollary}
Let $n\geq 2$ and let $b,f\in C^1(\mathbb{R})$ be such that they vanish at 0, with $f$ satisfying condition \hyperlink{H}{(H)}. Then, for all $\varphi\in L^2(\mathbb{R}^n)$ we have
\vspace{0.3cm}

\item[(i)] If $b$ satisfies hypothesis (i) or (ii) of Lemma \ref{T1} then
\begin{equation*}
\label{sugb2}
\n |x|^{\beta-1}|b'(t)|^{1/2}|D_x|^\beta e^{ib(t)|D_x|^2}\varphi\n_{L^2(\mathbb{R}_t\times\mathbb{R}^n_x)}
\leq C\n\varphi\n_{L^2(\mathbb{R}^n_x)},
\end{equation*}
for $1-n/2<\beta<1/2$, and
\begin{equation*}
\label{sugf2}
\n |x|^{\alpha-m/2}|b'(t)|^{1/2}|D_x|^{\alpha}e^{ib(t)|D_x|^{m}}\varphi\n_{L^2(\mathbb{R}_t\times\mathbb{R}^n_x)}\leq C\n\varphi\n_{L^2(\mathbb{R}^n_x)},
\end{equation*}
for $(m-n)/2<\alpha<(m-1)/2$. Finally, for $m>0$ and for any $\alpha,\beta>0$, we have
 $$\n |x|^{\beta-1}|b'(t)|^{1/2}|D_x|^\beta e^{ib(t)|D_x|^2}\varphi\n_{L^2(\mathbb{R}_t\times\mathbb{R}^n_x)}$$
 $$\leq C \n |x|^{\beta-1}|f'(t)|^{1/2}|D_x|^{m/2+\beta-1}e^{if(t)|D_x|^{m}}\varphi\n_{L^2(\mathbb{R}_t\times\mathbb{R}^n_x)}.$$
 
\item[(ii)] If $b$ and $c$ satisfy hypothesis (iii) of Lemma \ref{T1} then
\begin{equation*}
\label{sugb3}
\n |x|^{\beta-1}|c(t)b'(t)|^{1/2}|D_x|^\beta e^{ib(t)|D_x|^2}\varphi\n_{L^2(\mathbb{R}_t\times\mathbb{R}^n_x)}
\leq C\n\varphi\n_{L^2(\mathbb{R}^n_x)},
\end{equation*}
for $1-n/2<\beta<1/2$, and
\begin{equation*}
\label{sugf3}
\n |x|^{\alpha-m/2}|c(t)b'(t)|^{1/2}|D_x|^{\alpha}e^{ib(t)|D_x|^{m}}\varphi\n_{L^2(\mathbb{R}_t\times\mathbb{R}^n_x)}\leq C\n\varphi\n_{L^2(\mathbb{R}^n_x)},
\end{equation*}
for $(m-n)/2<\alpha<(m-1)/2$. Finally, for $m>0$ and for any $\alpha,\beta>0$, we have
 $$\n |x|^{\beta-1}|c(t)b'(t)|^{1/2}|D_x|^\beta e^{ib(t)|D_x|^2}\varphi\n_{L^2(\mathbb{R}_t\times\mathbb{R}^n_x)}$$
 $$\leq 
C \n |x|^{\beta-1}|f'(t)|^{1/2}|D_x|^{m/2+\beta-1}e^{if(t)|D_x|^{m}}\varphi\n_{L^2(\mathbb{R}_t\times\mathbb{R}^n_x)}.$$
\end{corollary}

\begin{corollary}
Let $b\in C^1(\mathbb{R})$ be vanishing at 0. Then, for all $\varphi\in L^2(\mathbb{R}^n)$, we have the following properties.
\vspace{0.3cm}

\item[(i)] Let $n\geq 1, m>0$ and $s>1/2$. If $b$ satisfies hypothesis (i) or (ii) of Lemma \ref{T1} then
\begin{equation*}
\n \la x_n\ra^{-s}|b'(t)|^{1/2}|D_n|^{(m-1)/2}e^{ib(t)|D_n|^{m}}\varphi\n_{L^2(\mathbb{R}_t\times\mathbb{R}^n_x)}\leq C \n \varphi\n_{L^2(\mathbb{R}^n)}.
\end{equation*}

Moreover, if $n\geq 2, m>0$ and $s>1/2$, then
\begin{equation*}
\n \la x_1\ra^{-s}|b'(t)|^{1/2}|D_n|^{(m-1)/2}e^{ib(t)D_1|D_n|^{m-1}}\varphi\n_{L^2(\mathbb{R}_t\times\mathbb{R}^n_x)}\leq C\n \varphi\n_{L^2(\mathbb{R}^n)}.
\end{equation*}
\item[(ii)] Let $n\geq 1, m>0$ and $s>1/2$. If $b$ and $c$ satisfy hypothesis (iii) of Lemma \ref{T1} then
\begin{equation*}
\n \la x_n\ra^{-s}|c(t)b'(t)|^{1/2}|D_n|^{(m-1)/2}e^{ib(t)|D_n|^{m}}\varphi\n_{L^2(\mathbb{R}_t\times\mathbb{R}^n_x)}\leq C \n \varphi\n_{L^2(\mathbb{R}^n)}.
\end{equation*}

Moreover, if $n\geq 2, m>0$ and $s>1/2$, then
\begin{equation*}
\n \la x_1\ra^{-s}|c(t)b'(t)|^{1/2}|D_n|^{(m-1)/2}e^{ib(t)D_1|D_n|^{m-1}}\varphi\n_{L^2(\mathbb{R}_t\times\mathbb{R}^n_x)}\leq C\n \varphi\n_{L^2(\mathbb{R}^n)}.
\end{equation*}
\end{corollary}

\begin{corollary}
Let $m>0$, $(m-n+1)/2<\alpha<(m-1)/2$ and let $b\in C^1(\mathbb{R})$ be vanishing at 0. Then, for all $\varphi\in L^2(\mathbb{R}^n)$, we have the following properties
\vspace{0.3cm}

\item[(i)] If $b$ satisfies hypothesis (i) or (ii) of Lemma \ref{T1} then
\begin{equation*}
\n  |x|^{\alpha-m/2}|b'(t)|^{1/2}|D_n|^{(m-1)/2}e^{ib(t)(|D_1|^m-|D'|^{m})}\varphi\n_{L^2(\mathbb{R}_t\times\mathbb{R}^n_x)}\leq C \n \varphi\n_{L^2(\mathbb{R}^n)},
\end{equation*}
where $D'=(D_2,...,D_n)$.
\item[(ii)] If $b$ and $c$ satisfy hypothesis (iii) of Lemma \ref{T1} then
\begin{equation*}
\n  |x|^{\alpha-m/2}|c(t)b'(t)|^{1/2}|D_n|^{(m-1)/2}e^{ib(t)(|D_1|^m-|D'|^{m})}\varphi\n_{L^2(\mathbb{R}_t\times\mathbb{R}^n_x)}\leq C \n \varphi\n_{L^2(\mathbb{R}^n)},
\end{equation*}
where $D'=(D_2,...,D_n)$.
\end{corollary}

\section{Weighted Strichartz estimate}\label{Strichartz}
In this section we shall prove that weighted Strichartz estimates hold for operators of the form $\LL_{b,\Delta}$ when $b$ is vanishing at 0 and satisfies the conditions considered before. The results will follow by using the well known results holding for classical Schr\"{o}dinger operators $\LL_{t,\Delta}$. 
\vspace{0.3cm}

\noindent{\bf Notations.} For brevity we shall denote $L^q_tL^p_x:=L^q_t(\mathbb{R};L^p_x(\mathbb{R}^n))$, and, sometimes, if not confusing, we shall use the same notation $L^q_tL^p_x:=L^q_t([0,T];L^p_x(\mathbb{R}^n))$ when considering a fixed time interval. Similarly for the mixed spaces $L^p_xL^q_t:=L^p_x(\mathbb{R}^n;L^q_t(\mathbb{R}))$ and  $L^p_xL^q_t:=L^p_x(\mathbb{R}^n;L^q_t([0,T]))$. Finally we shall denote $p'$ the conjugate exponent of $p$.

\begin{definition}[{\it Admissible pairs}]
Given $n\geq1$ we shall call a pair of exponents $(q,p)$ admissible if $2\leq q,p\leq \infty$, and 
$$\frac 2 q+\frac n p=\frac n 2, \quad \text{with}\quad (q,p,n)\neq (2,\infty,2).$$
\end{definition}

\begin{theorem}[{\it Global weighted Strichartz estimates}]\label{ws}
Let $b\in C^1(\mathbb{R})$ be vanishing at 0. Then, for any $(q,p)$,$(\tilde{q},\tilde{p})$ admissible pairs such that $2<q,\tilde q, p, \tilde p<\infty$, the following estimates hold
\vspace{0.3cm}

\item[(i)] if $b$ satisfies hypothesis (i) or (ii) of Lemma \ref{T1}, then we have
\vspace{0.3cm}

the weighted homogeneous Strichartz estimate 
\begin{equation}
\label{wh2}
\n |b'(t)|^{1/q}e^{ib(t)\Delta}\varphi\n_{L^q_tL^p_x}\leq C(n,q,p)\n\varphi\n_{L^2_x(\mathbb{R}^n)}, 
\end{equation}

the dual weighted homogeneous Strichartz estimate
\begin{equation}
\label{dwh2}
\n \int_\mathbb{R}|b'(s)|^{1/\tilde{q}}e^{-ib(s)\Delta}g(s)ds\n_{L^2_x}\leq C(n,\tilde{q},\tilde p)\n g\n_{L_t^{\tilde{q}'}L^{\tilde{p}'}_x},
\end{equation}

and the weighted inhomogeneous Strichartz estimate
\begin{equation}
\label{wi2}
\n |b'(t)|^{1/q}\int_\mathbb{R}|b'(s)|^{1/\tilde{q}}e^{i(b(t)-b(s))\Delta}g(s)ds\n_{L^q_tL^p_x}\leq C(n,q,p,\tilde{q},\tilde p)\n g\n_{L_t^{\tilde{q}'}L^{\tilde{p}'}_x}.
\end{equation}

\item[(ii)] if $b$ and $c$ satisfy hypothesis  (iii) of Lemma \ref{T1}, then we have
\vspace{0.3cm}

the weighted homogeneous Strichartz estimate 
\begin{equation}
\label{wh3}
\n |c(t)b'(t)|^{1/q}e^{ib(t)\Delta}\varphi\n_{L^q_tL^p_x}\leq C(n,q,p)\n\varphi\n_{L^2_x(\mathbb{R}^n)}, 
\end{equation}

the dual weighted homogeneous Strichartz estimate
\begin{equation}
\label{dwh3}
\n \int_\mathbb{R}|c(s)b'(s)|^{1/\tilde{q}}e^{-ib(s)\Delta}g(s)\n_{L^2_x}\leq C(n,q,p,\tilde{q},\tilde p)\n g\n_{L_t^{\tilde{q}'}L^{\tilde{p}'}_x},
\end{equation}

and the weighted inhomogeneous Strichartz estimate
\begin{equation}
\label{wi3}
\n |c(t)b'(t)|^{1/q}\int_\mathbb{R}|c(s)b'(s)|^{1/\tilde{q}}e^{i(b(t)-b(s))\Delta}g(s)ds\n_{L^q_tL^p_x}\leq C(n,q,p,\tilde{q},\tilde p)\n g\n_{L_t^{\tilde{q}'}L^{\tilde{p}'}_x}.
\end{equation}
\vspace{0.3cm}

\end{theorem}

\proof
The proof is a combination of standard results and duality arguments.
We first recall that, by Lemma \ref{T1}, we have
$$ \n b'(t)^{1/q}e^{ib(t)\Delta}\varphi\n_{L^q_tL^p_x}=\n b'(t)^{1/q}\n e^{ib(t)\Delta}\varphi\n_{L^p_x}\n_{L^q_t}\leq C \n e^{it\Delta}\varphi\n_{L^q_tL^p_x},$$
where $C$ depends on the properties of $b$ (see Lemma \ref{T1}).
Therefore we immediately get \eqref{wh2} from the standard case, that is, 
$$\n b'(t)^{1/q}e^{ib(t)\Delta}\varphi\n_{L^q_tL^p_x}
\leq C \n e^{it\Delta}\varphi\n_{L^q_tL^p_x}
\leq C(n,q,p)\n\varphi\n_{L^2_x(\mathbb{R}^n)},$$
where the last inequality corresponds to the classical homogeneous Strichartz estimate. Note that the application of the classical homogeneous Strichartz estimate, in particular, gives us the standard range for the admissible pairs.

As regards \eqref{dwh2}, it follows by duality and by \eqref{wh2} applied on the admissible pair $(\tilde{q},\tilde{p})$. Moreover, one can actually see (by duality) that  \eqref{dwh2} and \eqref{wh2} are equivalent.

As regards \eqref{wi2} it is given by duality, H\"{o}lder's inequality and \eqref{dwh2}. In fact, given $\psi\in L^{q'}_tL^{p'}_x$ we have

$$\int_{\mathbb{R}_t}\int_{\mathbb{R}_x^n} b'(t)^{1/q}\left(\int_\mathbb{R}b'(s)^{1/\tilde{q}}e^{i(b(t)-b(s))\Delta}g(s)ds\right) \overline{\psi(t,x)}dx\,\,dt$$
$$=\int_{\mathbb{R}_x^n}\left(\int_\mathbb{R}b'(s)^{1/\tilde{q}}e^{-ib(s)\Delta}g(s,x)ds\right)\left( \overline{\int_{\mathbb{R}_t} b'(t)^{1/q} e^{-ib(t)\Delta}\psi(t,x)dt}\right) dx$$
$$\leq \n  \int_\mathbb{R}b'(s)^{1/\tilde{q}}e^{-ib(s)\Delta}g(s,x)ds\n_{L^2_x}\n\int_{\mathbb{R}_t} b'(t)^{1/q} e^{-ib(t)\Delta}\psi(t,x)dt\n_{L^2_x}$$
$$\leq C(n,q,p,\tilde q, \tilde p) \n g\n_{L^{\tilde{q}'}_tL^{\tilde{p}'}}\n \psi\n_{L^{q'}_tL^{p'}_x},$$
which gives \eqref{wi2} (by duality). Moreover, one has, once again, equivalence between \eqref{wi2} and \eqref{dwh2}.
The other cases covered by \eqref{wh3},  \eqref{dwh3} and  \eqref{wi3} follow by using the same procedure.
\endproof

The advantage of the weighted Strichartz estimates above is that one can recover the standard range for the admissible exponents. Moreover, since the estimates are all equivalent, and since the endpoint case is solved in the standard Schr\"{o}dinger case (see \cite{KT}), which, in particular, immediately gives \eqref{wh2} and its generalizations, we can conclude that the estimates in Theorem \ref{ws}  hold true in the endpoint case $(q,p,n)=(2,\frac{2n}{n-2},n)$ as well.
We then have the following theorem.

\begin{theorem}[\it{Endpoint weighted Strichartz estimates}]
Let $b\in C^1(\mathbb{R})$ be such that it satisfies hypothesis (i), (ii), or (iii) of Lemma \ref{T1}. Then, if $b$ satisfies hypothesis (i) or (ii) of Lemma \ref{T1}, we have that \eqref{wh2}, \eqref{dwh2} and \eqref{wi2} are still satisfied for $(q,p,n)=(2,\frac{2n}{n-2},n)$ and/or $(\tilde{q},\tilde{p},n)=(2,\frac{2n}{n-2},n)$, with $n\geq 3$.
If $b$ satisfies hypothesis (iii) of Lemma \ref{T1} we have that \eqref{wh3}, \eqref{dwh3} and \eqref{wi3} are still satisfied for $(q,p,n)=(2,\frac{2n}{n-2},n)$ and/or $(\tilde{q},\tilde{p},n)=(2,\frac{2n}{n-2},n)$, with $n\geq 3$.

\end{theorem}

However, in order to use Strichartz-type estimates to prove the local well-posedness of the related semilinear Cauchy problem, we must be able to localize in time. In this sense the following estimates will be used.

\begin{theorem}[Local weighted Strichartz estimates]\label{swp}
Let $b\in C^1([0,T])$ be strictly monotone and vanishing at 0. Then, on denoting by $L^q_tL^p_x:=L^q([0,T];L^p(\mathbb{R}^n))$, we have that for any $(q,p)$ admissible pair such that $2< q,p< \infty$, the following estimates hold

\begin{equation}
\label{lwh2}
\n |b'(t)|^{1/q}e^{ib(t)\Delta}\varphi\n_{L^q_tL^p_x}\leq C(n,q,p)\n\varphi\n_{L^2_x(\mathbb{R}^n)}, 
\end{equation}

\begin{equation}
\label{lwh3}
\n e^{ib(t)\Delta}\varphi\n_{L^\infty_tL^2_x}\leq \n\varphi\n_{L^2_x(\mathbb{R}^n)}, 
\end{equation}

\begin{equation}
\label{lwi2}
\n |b'(t)|^{1/q}\int_0^t |b'(s)| e^{i(b(t)-b(s))\Delta}g(s)ds\n_{L^q_tL^p_x}\leq C(n,{q},p)\n |b'|^{1/{q'}}g\n_{L_t^{{q}'}L^{{p}'}_x},
\end{equation}
and
\begin{equation}
\label{lwi3}
\n \int_0^t |b'(s)| e^{i(b(t)-b(s))\Delta}g(s)ds\n_{L^\infty_tL^2_x}\leq C(n,{q},p)\n |b'|^{1/{q'}}g\n_{L_t^{{q}'}L^{{p}'}_x}.
\end{equation}

\end{theorem}

\proof
Estimate \eqref{lwh3}  trivially follows by the unitarity of $e^{ib(t)\Delta}$, whereas estimate \eqref{lwh2} follows by change of variables. In fact we have
$$\n |b'(t)|^{1/q}e^{ib(t)\Delta}\varphi\n_{L^q_tL^p_x}\leq \n e^{it\Delta}\varphi\n_{L^q_t([0, \tilde{T}];L^p_x)}
\lesssim_{n,q,p}\n \varphi\n_{L^2_x},$$
where in the last inequality we applied the standard homogeneous Strichartz estimate (notice that the first inequality is an equality if $b$ is strictly monotone).

As regards estimate \eqref{lwi2}, by using the change of variables $t'=b(t)$ and $s'=b(s)$, we have
$$\n |b'(t)|^{1/q}\int_0^t |b'(s)| e^{i(b(t)-b(s))\Delta}g(s)ds\n_{L^q_tL^p_x}\leq
\n\int_0^{b(t)} e^{i(t'-s')\Delta}\tilde{g}(s')ds'\n_{L^q_{t'}([0,\tilde{T}];L^p_x)}$$
$$= \n\int_0^{b(T)} e^{i(t'-s')\Delta}\chi(s')\tilde{g}(s')ds'\n_{L^q_{t'}([0,\tilde{T}];L^p_x)}$$
where $\tilde{g}=g\circ b^{-1}$, $\tilde{T}=b(T)$ and $\chi=1_{[0,b(t)]}$.
We then analyze the last quantity and, by using the properties of the Schr\"{o}dinger group $e^{it\Delta}$, we have
$$ \n\int_0^{b(T)} e^{i(t'-s')\Delta}\chi(s')\tilde{g}(s')ds'\n_{L^q_{t'}([0,\tilde{T}];L^p_x)}
\leq \n\int_0^{b(T)} \n e^{i(t'-s')\Delta}\chi(s)\tilde{g}(s')\n_{L^p_x}ds'\n_{L^q_{t'}([0,\tilde{T}])}$$
$$\leq \n\int_0^{b(T)}  \frac{1}{|t'-s'|^{n(1/2-1/p)}}\n \chi(s')\tilde{g}(s')\n_{L^{p'}_x}ds'\n_{L^q_{t'}([0,\tilde{T}])}$$
$$\underset{\text{H-L-S}} {\leq }C(n,q,p)\n \tilde{g}\n_{L^{q'}_{t'}([0,\tilde{T}];L^{p'}_x)}\underset{t=b^{-1}(t')}{\leq}C(n,q,p)\n |b'|^{1/{q'}}g\n_{L^{q'}_tL^{p'}_x},$$
where H-L-S stands for the application of the Hardy-Littlewood-Sobolev inequality, and we get \eqref{lwi2}.

We are now left with the proof of \eqref{lwi3}. We first consider the $L^2_x$-norm and have
$$\n \int_0^t |b'(s)| e^{i(b(t)-b(s))\Delta}g(s)ds\n_{L^2_x}^2=\n \int_0^t |b'(s)| e^{-ib(s))\Delta}g(s)ds\n_{L^2_x}^2$$
$$=\int_{\mathbb{R}^n} \left( \int_0^t |b'(s)| e^{-ib(s)\Delta}g(s)ds\right)\overline{\left( \int_0^t |b'(s')| e^{-ib(s')\Delta}g(s')ds'\right)}dx$$
$$= \int_0^t \int_{\mathbb{R}^n} |b'(s)|^{1/{q'}} g(s) \overline{\left(  |b'(s)|^{1/q}\int_0^t |b'(s')| e^{i(b(s)-b(s'))\Delta}g(s')ds'\right)}dx ds$$
$$\leq \int_0^t \n |b'(s)|^{1/{q'}} g(s) \n_{L^{p'}_x}\n |b'(s)|^{1/q}\int_0^t |b'(s')| e^{i(b(s)-b(s'))\Delta}g(s')ds'\n_{L^{p}_x} ds$$
$$\underset{(t<T)}{\leq}\n |b'|^{1/{q'}} g \n_{L^{q'}_sL^{p'}_x}\n |b'(s)|^{1/q}\int_0^t |b'(s')| e^{i(b(s)-b(s'))\Delta}g(s')ds'\n_{L^q_sL^{p}_x}$$
$$\underset{\text{by \eqref{lwi2}}}{\leq} C(n,q,p) \n |b'|^{1/{q'}} g \n_{L^{q'}_tL^{p'}_x}^2,$$
which, in particular, gives \eqref{lwi3}.
\endproof

We can also prove the previous result in a more general case, that is, when $b$ is not strictly monotone and has a finite number of critical points in the time interval $[0,T]$ (i.e. $b$ satisfies hypothesis of point (ii) in Lemma \ref{T1}).

\begin{theorem}\label{swp2}
Let $b\in C^1([0,T])$ be vanishing at 0 and such that $\sharp\{t\in[0,T]; b'(t)=0\}=k\geq 1$. Then, on denoting by $L^q_tL^p_x:=L^q([0,T];L^p(\mathbb{R}^n))$, we have that for any $(q,p)$ admissible pair such that $2< q,p< \infty$, the following estimates hold

\begin{equation}
\label{wh2.1}
\n |b'(t)|^{1/q}e^{ib(t)\Delta}\varphi\n_{L^q_tL^p_x}\leq C(n,q,p,k)\n\varphi\n_{L^2_x(\mathbb{R}^n)}, 
\end{equation}

\begin{equation}
\label{wh3.1}
\n e^{ib(t)\Delta}\varphi\n_{L^\infty_tL^2_x}\leq \n\varphi\n_{L^2_x(\mathbb{R}^n)}, 
\end{equation}

\begin{equation}
\label{wi2.1}
\n |b'(t)|^{1/q}\int_0^t |b'(s)| e^{i(b(t)-b(s))\Delta}g(s)ds\n_{L^q_tL^p_x}\leq C(n,{q},p,k)\n |b'|^{1/{q'}}g\n_{L_t^{{q}'}L^{{p}'}_x},
\end{equation}
and
\begin{equation}
\label{wi3.1}
\n \int_0^t |b'(s)| e^{i(b(t)-b(s))\Delta}g(s)ds\n_{L^\infty_tL^2_x}\leq C(n,{q},p,k)\n |b'|^{1/{q'}}g\n_{L_t^{{q}'}L^{{p}'}_x}.
\end{equation}

\end{theorem}

\proof
The proof uses the same arguments used in the proof of Theorem \ref{swp}. 
Estimate \eqref{wh3.1} is immediate. As regards \eqref{wh2.1}, we consider $0=T_0\leq T_1<T_2<...<T_k\leq T_{k+1}=T$ such that $b'(T_j)=0$ for $j=1,...k$, so that $b$ is strictly monotone on $[T_j,T_{j+1}]$, and we have
$$\n |b'(t)|^{1/q}e^{ib(t)\Delta}\varphi\n_{L^q_tL^p_x}=\left( \sum_{j=0}^k \n|b'(t)|^{1/q}e^{ib(t)\Delta}\varphi \n^q _{L^q([T_j,T_{j+1}];L^{p}_x)}\right)^{1/q}$$
$$\leq \sum_{j=0}^k \n|b'(t)|^{1/q}e^{ib(t)\Delta}\varphi \n _{L^q([T_j,T_{j+1}];L^{p}_x)}$$
$$\underset{b(t)=t'}{\leq} \sum_{j=0}^k \n e^{it\Delta}\varphi \n _{L^q([\tilde{T}_j,\tilde{T}_{j+1}];L^{p}_x)}$$
$$\leq  (k+1)C(n,q,p) \n \varphi\n_{L^2_x},$$
which proves the estimate.
We now prove \eqref{wi2.1} combining the previous procedure and the arguments in Theorem \ref{swp}.
By splitting the time interval again, we have
$$\n |b'(t)|^{1/q}\int_0^t |b'(s)| e^{i(b(t)-b(s))\Delta}g(s)ds\n_{L^q_tL^p_x}$$
$$\leq\sum_{j=0}^k \n |b'(t)|^{1/q}\int_0^t |b'(s)| e^{i(b(t)-b(s))\Delta}g(s)ds\n_{L^q_t([T_j,T_{j+1}];L^p_x)},$$
where, by the same procedure used in the proof of \eqref{lwi2}, each therm satisfies
$$\n |b'(t)|^{1/q}\int_0^t |b'(s)| e^{i(b(t)-b(s))\Delta}g(s)ds\n_{L^q_t([T_j,T_{j+1}];L^p_x)}$$
$$\lesssim \n |b'|^{1/{q'}}g\n_{L^{q'}_t([ T_j, T_{j+1}];L^{p'}_x)}\lesssim \n |b'|^{1/{q'}}g\n_{L_t^{{q}'}L^{{p}'}_x},$$
and we can conclude
$$\n |b'(t)|^{1/q}e^{ib(t)\Delta}\varphi\n_{L^q_tL^p_x}\leq (k+1)C(n,q,p) \n |b'|^{1/{q'}}g\n_{L_t^{{q}'}L^{{p}'}_x}.$$
Finally, estimate \eqref{wi3.1} follows by the same arguments used in the proof of \eqref{lwi3}.
\endproof

\section{Local well-posedness of the semilinear Cauchy problem}\label{LWP}

We can now apply the previous results to obtain the local well-posedness of the semilinear Cauchy problem
\begin{equation}
\label{NL1}
\left\{\begin{array}{l}
\partial_tu+ib'(t)\Delta u=\mu |b'(t)| |u|^{p-1}u,\\
u(0,x)=u_0(x),\\
\end{array}\right .
\end{equation}
where $b\in C^1([0,+\infty))$  vanishes at 0 and is either strictly monotone or has a finite number of critical points in any finite time interval $[0,T]$, for any $T<\infty$.

\begin{theorem}\label{TNL1}
Let  $1<p<\frac{4}{n}+1$ and $b\in C^1([0,+\infty))$ be vanishing at 0 and it is either strictly monotone or such that $\sharp\{t\in[0,\tilde{T}]; b'(t)=0\}$ is finite for any $\tilde{T}<\infty$. Then for all $u_0\in L^2(\mathbb{R}^n)$ there exists $T=T(\n u_0\n_2,n,\mu,p)>0$ such that there exists a unique solution $u$ of the IVP \eqref{NL1}
 in the time interval $[0,T]$ with
$$u\in C([0,T];L^2(\mathbb{R}^n))\bigcap L^q_t([0,T];L^{p+1}_x(\mathbb{R}^n))$$
and $q=\frac{4(p+1)}{n(p-1)}$. Moreover the map $u_0\mapsto u(\cdot,t)$, locally defined from $L^2(\mathbb{R}^n)$ to $C([0,T); L^2(\mathbb{R}^n))$, is continuous.
\end{theorem}

\proof[Proof of Theorem \ref{TNL1}]
We shall prove the result by using the standard fixed point argument.
We consider  the space 
$$X_T:=\{ u\in C([0,T];L^2(\mathbb{R}^n)\bigcap L^q_t([0,T];L^{p+1}_x(\mathbb{R}^n)); \n u\n_{X_T}<\infty\}$$
where
$$ \n u\n_{X_T}:=\n u\n_{L^\infty_t L^2_x}+\n |b'(t)|^{1/q}u\n_{L^{q}_tL^{p+1}_x},$$
and $L^q_tL^p_x:=L^q([0,T];L^p_x(\mathbb{R}^n))$, and prove that the map
$$\Phi(u):= e^{ib(t)\Delta}u_0+\mu \int_0^t |b'(s)|e^{i(b(t)-b(s))\Delta}|u|^{p-1}u \,ds$$
is a contraction.

Let $q=\frac{4(p+1)}{n(p-1)}$  so that $(q,p+1)$ is an admissible pair. Then, by Theorem \ref{swp} and Theorem \ref{swp2} (depending on the properties of $b=b(t)$), we have

$$\n |b'(t)|^{1/q}\Phi(u)\n_{L^{q}_tL^{p+1}_x}\leq C\n u_0\n_{L^2_x}+|\mu|\n|b'(t)|^{1/q} \int_0^t |b'(s)|e^{i(b(t)-b(s))\Delta}|u|^{p-1}u \,ds\n_{L^{q}_tL^{p+1}_x}$$
$$\leq C\n u_0\n_{L^2_x} +C|\mu|\n |b'(t)|^{1/q'}\n u\n^p_{L^{p+1}_x}\n_{L^{q'}_t} $$
$$= C\n u_0\n_{L^2_x} +C|\mu|\left\{\int_0^T |b'(t)| \n u\n^{pq'}_{L^{p+1}_x} dt\right\}^{1/{q'}} $$
$$ \leq C\n u_0\n_{L^2_x} +C|\mu| T^{1-n(p-1)/4}\left\{\int_0^T |b'(t)|^{\frac{q}{pq'}}\n u\n^{q}_{L^{p+1}_x} dt\right\}^{p/{q}}$$
$$=C\n u_0\n_{L^2_x} +C|\mu| T^{1-n(p-1)/4}\left\{\int_0^T |b'(t)|^{(1+\frac{q}{pq'}-1)}\n u\n^{q}_{L^{p+1}_x} dt\right\}^{p/{q}}$$
$$\leq C\n u_0\n_{L^2_x} +C|\mu| T^{1-n(p-1)/4}\left( \sup_{[0,T]}|b'(t)|^{(\frac{q}{pq'}-1)\cdot \frac p q}\right)\n |b'(t)|^{1/q} u\n^{p}_{L^q_tL^{p+1}_x}$$
$$\leq C\n u_0\n_{L^2_x} +C(T)|\mu| \n u\n_{X_T}^p,$$
where $1-n(p-1)/4>0$ and $\frac{q}{pq'}>1$ and $C(T):=T^{1-n(p-1)/4}\left( \sup_{[0,T]}|b'(t)|^{\frac{q}{pq'}-1}\right)\searrow 0$ as $T\searrow 0$.
Note also that, once again by Theorem \ref{swp} or Theorem \ref{swp2} and by the computation above we have
$$\n \Phi(u)\n_{L^\infty_tL^2_x}\underset{\text{by \eqref{wh3.1},\eqref{wi3.1}}}{\leq}  \n u_0\n_{L^2_x}+C|\mu|\n |b'|^{1/{q'}}|u|^p\n_{L^{q'}_tL^{(p+1)/p}_x},$$
$$\leq C\n u_0\n_{L^2_x} +C(T)|\mu| \n u\n_{X_T}^p,$$
where $C(T)$ is the same constant given before.
Putting the two estimates together we finally have that
$$\n\Phi(u)\n_{X_T}\leq C\n u_0\n_{L^2}+2C(T)|\mu| \n u\n_{X_T}^p,$$
with $C=C(n,q,p,k)$ and $k$ denoting the number of critical points of $b$ (see \eqref{wh2.1}). Therefore, given $B(R)\subset X_T$ a ball of radius $R$ in $X_T$, with $R=R(\n u_0\n_{L^2_x})$, we get,  for $T$ sufficiently small,
 that $\Phi$ sends $B(R)$ into itself. 

Now we prove that $\Phi$ is a contraction. Given $u,v$ solution of \eqref{NL1} with initial data $u_0$, and $u,v\in B(R)$, we have
$$\n \Phi(u)-\Phi(v)\n_{X_T}=|\mu|\n |b'(t)|^{1/q} \int_0^t |b'(s)|e^{i(b(t)-b(s))\Delta}(u|u|^{p-1}-v|v|^{p-1})ds\n_{L^\infty_tL^{2}_x}$$
$$+|\mu|\n |b'(t)|^{1/q} \int_0^t |b'(s)|e^{i(b(t)-b(s))\Delta}(u|u|^{p-1}-v|v|^{p-1})ds\n_{L^q_tL^{p+1}_x}$$
$$=|\mu|(I+II).$$
For the term $II$ we have
$$II\leq C \n |b'|^{1/{q'}}(|u|^{p-1}+|v|^{p-1})(u-v)\n_{L^{q'}_tL^{(p+1)/p}_x}$$
$$\leq C \left\{\int_0^T |b'(t)|(\n u\n_{L^{p+1}_x}^{p-1}+\n v\n_{L^{p+1}_x}^{p-1})^{q'}\n u-v\n^{q'} _{L^{p+1}_x} dt\right\}^{1/{q'}}$$
$$\leq C\left\{\int_0^T \left(|b'(t)|^{(1-q'/q)1/{q'}}(\n u\n_{L^{p+1}_x}^{p-1}+\n v\n_{L^{p+1}_x}^{p-1})\right)^{q'} \left(|b'(t)|^{1/q}\n u-v\n _{L^{p+1}_x} \right)^{q'}dt\right\}^{1/{q'}}$$
$$\leq C \left(\int_0^T\left(|b'(t)|^{(1-q'/q)1/{q'}}(\n u\n_{L^{p+1}_x}^{p-1}+\n v\n_{L^{p+1}_x}^{p-1})\right)^{q/q-2}dt\right)^{(q-2)/q}$$
$$\times \left(\int_0^T |b'(t)|^{1/q}\n u-v\n^q_{L^{p+1}_x} dt\right)^{1/{q}}$$
$$\leq C \left\{\left(\int_0^T\left(|b'(t)|^{(1/{q'}-1/q)1/(p-1)}\n u\n_{L^{p+1}_x}\right)^{q(p-1)/(q-2)}dt \right)^{(q-2)/q}\right.$$
$$+\left. \left(\int_0^T\left(|b'(t)|^{(1/{q'}-1/q)1/(p-1)}\n v\n_{L^{p+1}_x}\right)^{q(p-1)/(q-2)}dt\right)^{(q-2)/q}
 \right\}\n |b'|^{1/q}(u-v)\n_{L^q_tL^{p+1}_x}.$$
Since $(p-1)/(q-2)<1$, by Holder's inequality we have
$$ \left(\int_0^T\left(|b'(t)|^{(1/{q'}-1/q)1/(p-1)}\n u\n_{L^{p+1}_x}\right)^{q(p-1)/(q-2)}\right)^{(q-2)/q}$$
$$\leq T^{\alpha} \left(\int_0^T\left(|b'(t)|^{(1/{q'}-1/q)1/(p-1)}\n u\n_{L^{p+1}_x}\right)^{q}\right)^{(p-1)/q}$$
$$\leq T^{\alpha} \left(\int_0^T |b'(t)|^{(q-2)/(p-1)}\n u\n_{L^{p+1}_x}^{q}\right)^{(p-1)/q}$$
$$\leq T^{\alpha}\sup_{[0,T]}|b'|^{(q-2)/(p-1)-1}\left(\int_0^T |b'(t)|\n u\n_{L^{p+1}_x}^{q}\right)^{(p-1)/q}$$
$$=C(T) \n |b'|^{1/q}u\n_{L^q_tL^{p+1}_x}^{p-1}\leq  C(T) \n u\n^{p-1}_{X_T},$$
where $C(T):=T^{\alpha}\sup_{[0,T]}|b'|^{(q-2)/(p-1)-1}\searrow 0$ as $T\searrow 0$, and $\alpha=\left(\frac{p-1}{q-(p+1)}\right)^{-1}>0$.
Repeating the previous estimate for the term containing the $\n v\n_{L_x^{p+1}}$ we finally have
$$II\leq C(T)(\n u\n_{X_T}^{p-1}+\n v\n_{X_T}^{p-1})\n u-v\n_{X_T}.$$
By the previous computations it follows that the quantity $I$ can be estimated exactly as the quantity $II$, so we get
$$\n \Phi(u)-\Phi(v)\n_{X_T}\leq C(T)|\mu|(\n u\n_{X_T}^{p-1}+\n v\n_{X_T}^{p-1})\n u-v\n_{X_T}$$
$$\leq 2C(T) |\mu|R^{p-1}\n u-v\n_{X_T},$$
where, recall, $R=R(\n u_0\n_{L^2_x})$.
Finally, eventually by choosing $T$ smaller than before, we get that $\Phi$ is a contraction, so we have existence and uniqueness of the solution in $X_T$.

To see that the solution depends continuously from the initial data, we consider $u$ and $v$ solutions of \eqref{NL1} with initial data $u_0$ and $v_0$ respectively. Then 
$$u(t)-v(t)= e^{ib(t)\Delta}(u_0-v_0)+\mu  \int_0^t |b'(s)|e^{i(b(t)-b(s))\Delta}(u|u|^{p-1}-v|v|^{p-1})ds,$$
where, recall, $u\in B_{R_1}$ with $R_1=R_1(\n u_0\n_{L^2})=c_1\n u_0\n_{L^2_x}$, and $v\in B_{R_2}$ with $R_2=R_2(\n v_0\n_{L^2})=c_2\n v_0\n_{L^2_x}$ (with suitable constants $c_1$ and $c_2$).
Since, by the arguments above,
$$\n |b'|^{1/q}(u-v)\n_{L^q_tL^{p+1}_x} 
\leq C \n u_0-v_0\n_{L^2_x} +C|\mu| \n |b'|^{1/{q'}}(|u|^{p-1}+|v|^{p-1})(u-v)\n_{L^{q'}_tL^{(p+1)/p}_x}$$
$$\leq C \n u_0-v_0\n_{L^2_x} +C(T)|\mu|(\n u\n_{X_T}^{p-1}+\n v\n_{X_T}^{p-1})\n  |b'|^{1/q}(u-v)\n_{L^q_tL^{p+1}_x},$$
$$\leq C \n u_0-v_0\n_{L^2_x} +C(T)|\mu|(\n u_0\n_{L^2_x}+\n v_0\n_{L^2_x})^{p-1}\n  |b'|^{1/q}(u-v)\n_{L^q_tL^{p+1}_x},$$
we have
$$\left(1-C(T)|\mu|(\n u_0\n_{L^2_x}+\n v_0\n_{L^2_x})^{p-1}\right )\n |b'|^{1/q}(u-v)\n_{L^q_tL^{p+1}_x}\leq C \n u_0-v_0\n_{L^2_x}.$$
Then, if $\n u_0-v_0\n_{L^2_x}$ is small enough, we get
$$1-C(T)|\mu|(\n u_0\n_{L^2_x}+\n v_0\n_{L^2_x})^{p-1}\geq 1-C(T)|\mu|\n u_0-v_0\n_{L^2_x}^{p-1}\geq c_0>0$$
which gives 
$$ \n |b'|^{1/q}(u-v)\n_{L^q_tL^{p+1}_x}\leq c \n u_0-v_0\n_{L^2_x},$$
where $c>0$ is a new suitable constant.
Finally, since
$$\n u-v\n_{L^\infty_t L^2_x}\leq C \n u_0-v_0\n_{L^2_x} +C(T)|\mu|\n |b'|^{1/{q'}}(|u|^{p-1}+|v|^{p-1})(u-v)\n_{L^{q'}_tL^{(p+1)/p}_x}$$
by the previous computations we have
$$\n u-v\n_{L^\infty_t L^2_x}\leq \tilde C \n u_0-v_0\n_{L^2_x},$$
which concludes the proof.
 \endproof

\begin{remark}
The previous result shows that, if either the initial data or the time $T$ is small enough, then we have local well-posedness of the IVP \eqref{NL1} even in presence of a function $b'$ vanishing at a finite number of times in $[0,T]$. In particular, if $b'$ vanishes at more that one point and we consider a very small initial data,  then the solution very likely exists in a time interval $[0,T]$ containing more than one degenerate point. Conversely, if $b'$ vanishes at more than one point and  $T$ is very small, then the picture may coincide with the one we observe when $b$ is strictly monotone, since we may cross no more than one time-degeneracy.
\end{remark}

Let us now show a few examples of operators that fall in our case. 

\begin{example}
The most natural example we can give is the one generalizing the standard Schr\"{o}dinger operator, that is, when $b(t)=t^{\alpha+1}/(\alpha+1)$ with $\alpha\geq 0$. In fact, for $\alpha=0$, we recover exactly the classical case.
For any $\alpha\geq 0$, $b$ is strictly monotone on each finite time interval $[0,T]$, therefore the IVP \eqref{NL1} associated with the operator
$$\LL_{b,\Delta}=\LL_{\frac{t^{\alpha+1}}{\alpha+1},\Delta}=\partial_t+ it^\alpha \Delta,\quad $$
is locally solvable by means of Theorem \ref{TNL1}. 

It is worth to mention that operators of this form have been analyzed by Cicognani and Reissig in \cite{CR} who proved the local well posedness of the associated linear Cauchy problem both in Sobolev and Gevrey spaces. Moreover, the first author and Gigliola Staffilani  proved in \cite{FS} the local smoothing effect for operators of this form (also applicable to more general $b(t)$) and proved the local well-posedness of the associated nonlinear Cauchy problem both with polynomial and with derivative nonlinearities.
\end{example}

\begin{example}
An other example that falls in the case when $b$ is strictly monotone in $[0,T]$, for any $0<T<\infty$, is when $b(t)=e^{t}-t-1$. Once again we have that the IVP \eqref{NL1}, with
$$\LL_{b,\Delta}=\LL_{e^{t}-t-1,\Delta}=\partial_t+i(e^{t}-1)\Delta,$$
is locally solvable. Note also that the operator is degenerate at $t=0$.
\end{example}

\begin{example}
In this example we consider $b(t)=\sin(t)$. This function satisfies our hypotheses since it vanishes at 0 and 
$\sharp\{t\in[0,T]; b'(t)=0\}=k\geq 1$ for any $0<T<\infty$. Then the IVP \eqref{NL1} for the operator
$$\LL_{b,\Delta}=\LL_{\sin(t),\Delta}:=\partial_tu+i\cos(t)\Delta $$
is locally solvable by Theorem \ref{TNL1}.
Similarly, one can consider other trigonometric functions, as, for instance, $b(t)=\cos(t)-1$, $b(t)=\sin(t)\cos(t)$ and more complicated ones. 
\end{example}

\end{document}